\documentclass[leqno,twoside]{amsart}

\usepackage{amssymb}
\usepackage{mathrsfs}
\usepackage{enumerate}
\usepackage{amsmath,amscd}
\usepackage{latexsym,esint}
\usepackage{color}
\usepackage{lscape}
\usepackage{rotating}

\usepackage{tikz}
\usetikzlibrary{arrows}
\usepackage[all]{xy}
\xyoption{2cell}
\xyoption{curve}
\UseTwocells
\input{diagxy}

\theoremstyle{theorem}
\newtheorem{theorem}{Theorem}[subsection]

\newtheorem{lemma}[theorem]{Lemma}
\newtheorem{proposition}[subsection]{Proposition}
\newtheorem{remark}[theorem]{Remark}

\theoremstyle{definition}
\newtheorem{definition}[subsection]{Definition}

\numberwithin{equation}{section}
\numberwithin{figure}{section}

\newcommand{\n}{[n]}

\newcommand{\tto}{\Rightarrow}
\newcommand{\iiota}{\overline{\iota}}
\newcommand{\nin}{\not\in}

\newcommand{\aaa}{\mathsf{A}}
\newcommand{\bbb}{\mathsf{B}}
\newcommand{\ccc}{\mathsf{C}}
\newcommand{\ddd}{\mathsf{D}}

\newcommand{\ctn}{\mathbb{C}}
\newcommand{\hcn}{N_\Delta(N_2\mathsf{Cat})}

\newcommand{\minn}[1]{p}
\newcommand{\maxx}[1]{q}
\newcommand{\spi}[1]{\mathrm{sp}(#1)}
\newcommand{\skn}{N_s(\mathsf{Cat})}
\newcommand{\ds}{X}
\newcommand{\dom}{\mathsf{dom}}
\newcommand{\cod}{\mathsf{cod}}

\begin{document}

\title[]{The Catalan simplicial set and uniform classification of monoidal-type categories}
\author{Aaron Greenspan} 
\maketitle

\begin{abstract}
  
Many monoidal-type objects are known to be classified by maps from the Catalan simplicial set $\ctn$ to various nerves of categories and higher categories.  There are, for example, three different nerves of the 2-category of categories $\mathsf{Cat}$ with the property that maps from $\ctn$ into them are in one to one correspondence with strict monoidal, monoidal, and skew monoidal categories.  In this last case, the classification result can be understood as justifying the definition of skew monoidal categories intrinsically in terms of the structure of $\ctn$.

In this paper, we consider a single nerve -- the homotopy coherent nerve of $\mathsf{Cat}$ viewed as a one object simplicially enriched category -- and verify that maps from $\ctn$ to it classify the three monoidal-type categories mentioned above, as well as lax monoidal categories and $\Sigma$-monoidal categories.  Furthermore, we identify a new monoidal-type category corresponding to the most general such maps and present them in the usual way: a collection of data subject to explicit coherence requirements.   These results suggest that the space of maps from $\ctn$ to the coherent nerve of $\mathsf{Cat}$ provides a uniform framework for the definition and study of many species of monoidal-type category, old or new.

\end{abstract}

\medskip 

\tableofcontents

\medskip
  
\section{Introduction}

The presence of natural transformations in the $2$-category of categories $\mathsf{Cat}$ gives a rich vocabulary for describing the ways in which a ``tensor product" functor $\otimes:A\times A\to A$ for a category $A$ may be associative and unital.  We have firstly the \textit{strict} monoidal categories, in which $\otimes$ is \textit{strictly} associative and unital;  weakening the definition slightly gives monoidal categories in the sense of \cite{sixth}, in which $\otimes$ is associative and unital only up to coherent natural isomorphism;   and one possible further weakening yields \textit{skew} monoidal categories -- introduced by Szlach\'anyi in \cite{first} -- in which $\otimes$ is associative and unital up to coherent, but not necessarily invertible, natural transformation.   

\medskip

In the latter two cases, it is worth saying something about our use of the word `coherent'.  For a monoidal category, the coherence of the associativity and unitality natural isomorphisms means that they render commutative the so-called pentagon and triangle diagrams.  According to a well known coherence theorem of Maclane's \cite{second}, the commutativity of these diagrams implies the commutativity of any well-formed diagram built out of these natural isomorphisms.  For skew monoidal categories, the coherence conditions require the commutativity of not two but five diagrams; and by contrast to the monoidal case, the commuativity of these diagrams does not imply the commutativity of all diagrams built from the natural transformations.   So where do these five diagrams come from? Why are they the `right' notion of coherence for skew monoidal categories?

\medskip

One possible answer to this question is given in \textit{The Catalan simplicial set} \cite{third} ; a key insight is that the theory of \textit{simplicial sets} provides a uniform framework for describing the data and coherence of monoidal-type categories.  For example, the data and coherence of a monoidal category consists in: a category, a pair of functors, three natural isomorphisms, and two commutative diagrams thereof.  As we go down this list we see that functors mediate between categories, natural transformations mediate between functors, and that commutative diagrams mediate between transformations.  Such structure suggests the simplices and face maps of a simplicial set, and indeed this can be formalized in terms of the \textit{psuedo nerve} $N_p(\mathsf{Cat})$ of the monoidal 2-category $\mathsf{Cat}$.  

\medskip

Explicitly, $N_{p}(\mathsf{Cat})$ is a simplicial set with a single $0$-simplex and (small) categories for $1$-simplices.  The $2$-simplices are binary functors $T:A\times B \to C$ where $A, B,$ and $C$ are the functor's three faces. The $3$-simplices are natural isomorphisms filling in squares of four such functors (the four faces) and higher simplices are commutative diagrams of such natural transformations.  All the data and coherence of a monoidal category live within $N_{p}(\mathsf{Cat})$: Specifically, a $1$-simplex, a pair of $2$-simplices, three $3$-simplices, and a pair of $4$-simplices, themselves suitably related by face maps.  By changing the strictness of the $3$-simplices, we obtain nerves $N_{\otimes}(\mathsf{Cat})$ and  $\skn$ suited to the description of strict and skew monoidal categories respectively.  For example, the five coherence diagrams in a skew monoidal category can be understood as five $4$-simplices in $\skn$.  

\medskip

We may wonder whether strict, skew, or plain monoidal categories are shadows of a fixed simplicial set, just as paths in a topological space are shadows of the interval.   It is shown in \cite{third} that this is indeed the case.  It defines the \textit{Catalan simplicial set} $\ctn$, recalled in section $2$ below, and shows that simplicial maps from $\ctn$ to $N_{p}(\mathsf{Cat})$, respectively $N_{\otimes}(\mathsf{Cat})$ and  $\skn$,  are in one to one correspondence with monoidal, respectively strict and skew monoidal, categories.   Now the structure of $\ctn$ itself provides an \textit{a priori} justification for the five coherence diagrams for skew monoidal categories: they are just what is required for a map $\ctn \to \skn$.   

\medskip

Many other kinds of `monoidal object' contained in a (higher-)categorical structure can be classified by maps out of $\ctn$ into suitably chosen nerves.  As shown in \cite{third} and \cite{fourth}, we may classify each of the following structures in this way: 

\medskip

\begin{enumerate}
\item Monoids in a fixed monoidal category (including the case of strict monoidal categories -- which are monoids in $\mathsf{Cat}$).
\item Monoidal categories.
\item Skew monoidal categories.
\item Monads in a bicategory.
\item Monoidales in a monoidal bicategory.
\item Skew monoidales in a monoidal bicategory.
\item Monoidal bicategories.
\item Skew monoidal bicategories.\footnote{The result about skew monoidal bicategories holds only after a mild restriction is placed on the maps out of $\ctn$ ; see \cite{fourth} Section 5.2.}
\end{enumerate}

\medskip

As with skew monoidal categories, these classifying results simultaneously justify otherwise complex and ad-hoc definitions.  It also suggests that new kinds of monoidal object may be defined directly in terms of maps from $\ctn$ into any reasonably defined nerve. 
\medskip

The story of this paper begins with a missing item from the above list: is it possible to exhibit monoidal $(\infty,1)$-categories in the sense of \cite{seventh} as maps out of $\ctn$?  The idea is as follows.  There is a quasi-category $\mathsf{qCat}$ whose 0-cells are small quasi-categories, whose 1-cells are quasi-functors, and whose higher cells are suitable quasi-invertible transformations.  If we set aside problems of strict associativity, $\mathsf{qCat}$ becomes a simplicial monoid under binary product of quasi-categories, and so can be seen as a one-object simplicially enriched category.  By taking its homotopy coherent nerve -- originally introduced in \cite{eighth} -- we obtain a simplicial set $N_{\Delta}(\mathsf{qCat})$ that should contain all the data required to define a monoidal $(\infty,1)$-category, so that these should be definable in terms of maps $\ctn \to N_{\Delta}(\mathsf{qCat})$.  There is, however, an analogous construction in ordinary category theory which deserves attention first, and which is the main subject of our present investigation.

\medskip

Like $\mathsf{qCat}$, the $2$-category $\mathsf{Cat}$ of categories, functors, and natural transformations can be viewed under the nerve construction as a simplicial set.  If we once again ignore problems of strict associativity, this simplicial set is a simplicial monoid under product of categories, and hence as a one-object simplicially enriched category, which we will refer to as $N_2(\mathsf{Cat})$.   Applying the homotopy coherent nerve yields a simplical set $\hcn$, and we can now consider maps $\ctn \to \hcn$.   Far from being a simple warm up for the higher categorical case, these maps hold great interest in their own right: in some sense, they unify an extensive array of monoidal-type categories. 

\medskip

As strict monoidal categories include into monoidal categories, and monoidal categories include into skew monoidal categories \cite{first}, so the three corresponding nerves admit parallel inclusions: $N_{\otimes}(\mathsf{Cat}) \subseteq N_{p}(\mathsf{Cat}) \subseteq\skn$.  Thus each strict monoidal, monoidal, or skew monoidal category can be understood as an element of the set $\mathsf{sSet}(\ctn,\skn)$.  On the contrary, not \textit{every} sort of monoidal-type category can be understood in this way; there is another common variant of monoidal structure which evades classification by $\skn$ and $\ctn$.  The definition of a \textit{lax} monoidal structure on a category $A$ approaches the idea of weakened associativity by introducing $n$-ary operations $\otimes^n:A^n \to A$ for each $n\ge 0$, which we think of as ``parenthesis free'' $n$-ary multiplications.  As the associativity of a binary operation is given by relationships between higher arity composites of itself, the functors $\otimes^n$ are used to mediate between possible such composites.  As before, strict monoidal implies monoidal implies lax monoidal, but there are lax monoidal categories which are not skew and vice versa.  Consequently, there are lax monoidal categories which cannot be understood as maps $\ctn \to \skn$.   This is where  $\hcn$ comes into play.

\medskip

In Section $4$ of this paper, we assign a simplicial map $\ctn \to \hcn$ to each lax monoidal category in such a way that the original lax monoidal category can be completely reconstructed from the assignment.  This gives an injective map taking lax monoidal categories into $\mathsf{sSet}(\ctn,\hcn)$.  Though we do not explicitly characterize the image of our assignment, it shows that  $\mathsf{sSet}(\ctn,\hcn)$ can classify lax monoidal categories in the spirit of the results of \cite{third} and \cite{fourth}.  In Section $5$ we show that the nerve $\skn$ classifying skew monoidal categories also includes into $\hcn$.  This shows that each kind of monoidal category mentioned thus far can be identified with a map $\ctn \to \hcn$.  In Section $6$ we examine other kinds of monoidal categories arising as maps in $\mathsf{sSet}(\ctn,\hcn)$ and among them find $\Sigma$-\textit{monoidal} categories as in \cite{fifth}.  In this sense, the set $\mathsf{sSet}(\ctn,\hcn)$ is rich enough to classify an extensive spectrum of monoidal-type structures in category theory. Furthermore, in sections $3$ and $6$ we see that there are maps $\ctn \to \hcn$ corresponding to yet a weaker type of monoidal-type category not yet defined in the literature. 

\medskip

Future work will include an examination of maps from $\ctn$ into the homotopy coherent nerve of various higher and enriched categories in the hopes of capturing a full range of possible monoidal objects in many more contexts.  This of course includes the original motivation where the category in question is that of quasi-categories. 

\medskip

\noindent\textbf{Acknowledgements:}  I would like to thank Richard Garner for drawing my attention to this topic and for his illuminating comments, suggestions, and advice throughout the investigation. 

\section{Definitions and notation: $\ctn$ and $\hcn$}

In this section we recall important definitions and introduce some helpful notations and shorthands.  We write $I$ for the terminal category, and $\Delta$ for the \textit{simplicial category} with finite ordered sets $[n]:= \lbrace 0<1<...<n \rbrace$ for objects, and order preserving maps.  These morphisms are generated by the \textit{coface} maps $\delta_i : [n-1] \to [n]$ and \textit{codegeneracy} maps $\sigma_i : [n+1] \to [n]$ given by:

\begin{displaymath}
    \delta_i(p):= \left\{
     \begin{array}{lr}
	    p & \text{if } p<i
		\vspace{.2cm}
		\\
		p+1 & \text{if } p\ge i 
     \end{array}
   \right.
   \hspace{1cm} \sigma_i(p):= \left\{
     \begin{array}{lr}
	    p & \text{if } p\le i
		\vspace{.2cm}
		\\
		p-1 & \text{if } p> i 
     \end{array}
   \right.
\end{displaymath} 

\medskip

\noindent We write $\mathsf{sSet}$ for the functor category $\mathsf{Fun}(\Delta^{op},\mathsf{Set})$ and refer to its objects as \textit{simplicial sets}.  A simplicial set $X:\Delta^{op} \to \mathsf{Set}$ thus gives rise to an indexed sequence of sets $X_n:=X[n]$,  \textit{face} maps $d_i = X\delta_i$, and \textit{degeneracy} maps $s_i = X\sigma_i$. As $X$ is contravariant, the face and degeneracy maps have now reversed direction: $d_i:X_n \to X_{n-1}$ and $s_i:X_n \to X_{n+1}$.  By slight abuse we may write `face' (`coface', `degeneracy', `codegeneracy') map to refer to a composite of face (coface, degeneracy, codegeneracy) maps.  This should always be made clear by the context.  At times we may simply use $d$ or $s$ to refer to a non-specific face or degeneracy map, and will always take $\delta$ and $\sigma$ to be the corresponding coface and codegeneracy map.

\medskip

An element $x\in X_n$ is an \textit{$n$-simplex}, and is referred to as \textit{degenerate} whenever $x$ is in the image of a degeneracy map.  When $m \le n$,  \textit{$x\in X_n$ is the degeneracy of $y\in X_m$} when $x$ is the image of $y$ under an $n-m$ fold composition of degeneracy maps $X_m \to X_n$.  We refer to the $(n-1)$-simplex $d_i(x)$ as the $i^{th}$\textit{face} of $x$.  More generally, we can identify arbitrary subsets $\ccc = \lbrace c_0<c_1<...<c_m \rbrace \subseteq [n]$ with the morphism $\delta_{\ccc}: [m] \to [n]$ of $\Delta$ which sends $i$ to $c_i$. For $x \in X_n$, denote $x_{\ccc} := X(\delta_{\ccc})(x) = d_{\ccc}x  \in X_m$. We think of and refer to this $m$-simplex as the \textit{$\ccc^{th}$-face} of $x$.\footnote{ With this notation, we might equally well write the face $d_i(x)$ as $x_{[n]-i}$.}   For two element subsets $\lbrace p<q \rbrace \subseteq [n]$, we drop the brackets and write simply $x_{p,q}$ or $x_{pq}$ for the $\lbrace p,q \rbrace^{th}$ face of $x$. In this notation, the \textit{spine} of the simplex $x$ is the collection of successive $1$-faces $\spi{x} : =\lbrace x_{01},x_{12},...,x_{n-1,n}\rbrace$.  It will be convenient to introduce the set $\ccc^-:= \ccc - \lbrace \max{\ccc} \rbrace$ and successor function $s:\ccc^- \to \ccc$ taking $c$ to the next greatest element of $\ccc$.  In this notation, $\spi{x_{\ccc}} = \lbrace x_{c,sc} ~|~ c \in \ccc^- \rbrace $.

\medskip

For an $n$-simplex $x \in X_n$, the \textit{boundary} of $x$ consists in the collection of its faces $\lbrace d_0(x),d_1(x),...,d_n(x)\rbrace$. Commutativity relations among coface maps $\delta_i$ give rise to relations amongst the faces of the boundary : $d_j(d_i(x))= d_j(d_{i+1}(x))$ for all $0\le i \le j <n$.  An \textit{$n$-boundary} in $X$ is a collection of ($n-1$)-simplices $\lbrace x_0,...,x_n \rbrace$ satisfying the same relationship: $d_j(x_i)=d_j(x_{i+1})$ for all $0\le i \le j <n$.  An $n$-boundary may in general be the boundary of one, many, or no $n$-simplices.  A simplicial set is called \textit{$r$-coskeletal} when for each $n>r$ and each $n$-boundary $\lbrace x_0,...,x_n \rbrace$, there is a unique $n$-simplex $x$ with $d_i(x) = x_i$ for $0\le i \le n$.  This establishes a bijection between $n$-boundaries, which are collections of $(n-1)$-simplices, and $n$-simplices.  Thus a definition of an $r$-coskeletal simplicial set need only specify the simplices up to dimension $r$, as all higher dimensional simplices are then determined by these bijections.  Similarly, a map from a simplicial set $Y$ into an $r$-coskeletal simplicial set $X$ is given by a map from the $r$-th truncation of $Y$ to $X$.  In otherwords, such a map is determined by where it sends simplices of dimension $\le r$.

\medskip

\begin{definition}
The \textit{Catalan simplicial set} $\ctn$ is the simplicial set defined by the following data:

\begin{enumerate}[\textbf{$\bullet$}]

\item A unique $0$-simplex: $*$.
\item Two $1$-simplices: $0=s_0(*)$, and $1$.
\item Five $2$-simplices:

\begin{center}

\begin{tikzpicture}[node distance=1cm]

  \node (A1) {$ *$} ;
  \node (A2) [right of=A1] {$ $};
  \node (B2) [below of=A2] {$*$};
  \node (B1) [below of=A1] {$ $};
  \node (B0) [left of=B1] {$*$};
  
  \draw[->, above right] (A1) to node {$ 0$} (B2);  
  \draw[->, above left] (B0) to node {$0$} (A1);
  \draw[->, below ] (B0) to node {$ 0 $} (B2);
    \draw[->, above ] (B0) to node {$ s_0(0) $} (B2);

\end{tikzpicture} , \begin{tikzpicture}[node distance=1cm]

  \node (A1) {$ *$} ;
  \node (A2) [right of=A1] {$ $};
  \node (B2) [below of=A2] {$*$};
  \node (B1) [below of=A1] {$ $};
  \node (B0) [left of=B1] {$*$};
  
  \draw[->, above right] (A1) to node {$ 1$} (B2);  
  \draw[->, above left] (B0) to node {$0$} (A1);
  \draw[->, below ] (B0) to node {$ 1 $} (B2);
  \draw[->, above ] (B0) to node {$ s_0(1) $} (B2);
\end{tikzpicture}, \begin{tikzpicture}[node distance=1cm]

  \node (A1) {$ *$} ;
  \node (A2) [right of=A1] {$ $};
  \node (B2) [below of=A2] {$*$};
  \node (B1) [below of=A1] {$ $};
  \node (B0) [left of=B1] {$*$};
  
  \draw[->, above right] (A1) to node {$ 0$} (B2);  
  \draw[->, above left] (B0) to node {$1$} (A1);
  \draw[->, below ] (B0) to node {$ 1 $} (B2);
  \draw[->, above ] (B0) to node {$ s_1(1) $} (B2);
\end{tikzpicture},\begin{tikzpicture}[node distance=1cm]

  \node (A1) {$ *$} ;
  \node (A2) [right of=A1] {$ $};
  \node (B2) [below of=A2] {$*$};
  \node (B1) [below of=A1] {$ $};
  \node (B0) [left of=B1] {$*$};
  
  \draw[->, above right] (A1) to node {$ 0$} (B2);  
  \draw[->, above left] (B0) to node {$0$} (A1);
  \draw[->, above ] (B0) to node {$ u_{} $} (B2);
  \draw[->, below ] (B0) to node {$ 1 $} (B2);
\end{tikzpicture} ,\begin{tikzpicture}[node distance=1cm]

  \node (A1) {$ *$} ;
  \node (A2) [right of=A1] {$ $};
  \node (B2) [below of=A2] {$*$};
  \node (B1) [below of=A1] {$ $};
  \node (B0) [left of=B1] {$*$};
  
  \draw[->, above right] (A1) to node {$ 1$} (B2);  
  \draw[->, above left] (B0) to node {$1$} (A1);
  \draw[->, below ] (B0) to node {$ 1 $} (B2);
  \draw[->, above ] (B0) to node {$ m_{} $} (B2);
\end{tikzpicture}
\end{center}

Which we may also write as:
\[s_0(0):0\vee 0 \to 0 ~,~ s_0(1):0 \vee 1 \to 1~,~ s_1(1):1 \vee 0 \to 1~, ~u: 0\vee 0 \to 1~,~ m:1 \vee 1 \to 1 \]

\item Higher-dimensional simplices are determined by $2$-coskeletality.

\end{enumerate}

\end{definition}

\medskip

Importantly for the proofs to follow, the $k$-simplices of $\ctn$ for $k\ge 3$ are determined by their $2$-simplex faces by coskeletality, and the $2$-simplices are determined by their $1$-simplex faces by definition.  As a result, all simplices of $\ctn$ are determined by their collection of $1$-simplex faces.  Thus we may identify $x$ with the set $\lbrace x_{p,q} ~|~ 0 \le p < q \le n \rbrace$ and $x_\ccc$ with $\lbrace x_{c,c'} ~|~ c < c' \in \ccc \rbrace$.

\medskip

The intuition behind the remarkable classifying properties of $\ctn$ is revealed by thinking of $2$-simplices in $\ctn$ as maps and higher simplices as diagrams of such maps. If we think of the non-degenerate $2$-simplex $m : 1 \vee 1 \to 1$ as a `monoidal product' of some kind, then the $3$-simplex filling the boundary $\lbrace m, m, m, m \rbrace$ can be thought of as encoding some kind of associativity of $m$.  

\medskip
\begin{center}
\begin{tikzpicture}[node distance=2cm]

  \node (A1) {$1 \vee 1 \vee 1$} ;
  \node (A2) [right of=A1] {$1\vee 1$};
  \node (B2) [below of=A2] {$1$};
  \node (B1) [below of=A1] {$1\vee 1$};
  \node (B0) [left of=B1] {$ $};
  
  \draw[->, above] (A1) to node {$ 1\vee m$} (A2);  
  \draw[->, left] (A1) to node {$m \vee 1$} (B1);
  \draw[->, below ] (B1) to node {$ m $} (B2);
    \draw[->, right ] (A2) to node {$ m $} (B2);

\end{tikzpicture}
\end{center}
\medskip

These ideas are made precise in a variety of contexts in \cite{third}.  Among their many results are those concerned with the skew monoidal categories of the introduction.

\medskip

\begin{definition}
A \textit{skew monoidal category} $(A, \otimes, \otimes^0, \alpha, \lambda, \rho)$ consists in:

\begin{enumerate} [\textbf{$\bullet$}]

\item A category $A$.
\item A functor $\otimes: A \times A \to A$.
\item A functor $\otimes^0:I \to A$. (i.e an object of $A$)
\item A natural tranformation:
\[ \alpha: \otimes \circ (\otimes \times 1_A) \tto \otimes \circ (1_A \times \otimes)\]
\item Natural transformations:
\[ \lambda: \otimes \circ (\otimes^0 \times 1_A) \tto 1_A \hspace{.25cm},\hspace{.25cm} \rho: 1_A \tto \otimes \circ (1_A \times \otimes^0)\] 

\end{enumerate}

\medskip

\noindent  These natural transformations must additionally commute in five diagrams.
\begin{description}
{ 
\item[Associativity]\hfill\\ The following pentagon commutes where each edge involves a single application of $\alpha$:

\begin{center}\begin{tikzpicture}[node distance=1cm]

  \node (A1) {$ $} ;
  \node (A2) [right of=A1] {$ $};
  \node (A3) [right of=A2] {$\otimes \circ (\otimes \times 1) \circ (\otimes \times 1 \times 1) $};
  
  \node (B1) [below of=A1] {$\otimes \circ (\otimes \times 1) \circ (1 \times \otimes \times 1)$};
  \node (B2) [right of=B1] {$ $};
  \node (B3) [right of=B2] {$ $};
  
  \node (C1) [below of=B1] {$ $};
  \node (C2) [right of=C1] {$ $};
  \node (C3) [right of=C2] {$ $};
  \node (C4) [right of=C3] {$ $};
  \node (C5) [right of=C4] {$\otimes \circ (\otimes \times \otimes) $};
  
  \node (D1) [below of=C1] {$ \otimes \circ (1 \times \otimes) \circ (1 \times \otimes \times 1)$};
  \node (D2) [right of=D1] {$ $};
  \node (D3) [right of=D2] {$ $};
  
  \node (E1) [below of=D1] {$ $};
  \node (E2) [right of=E1] {$ $};
  \node (E3) [right of=E2] {$ \otimes \circ (1 \times \otimes) \circ (1 \times 1 \times \otimes)$};

  \draw[->,  right] (A3) to node {$ $} (B1);
  \draw[->,  right] (A3) to node {$ $} (C5);
  \draw[->,  below right] (B1) to node {$ $} (D1);
  \draw[->, above left] (C5) to node {$ $} (E3);
  \draw[->, below left] (D1) to node {$ $} (E3);
\end{tikzpicture}\end{center}

\item[Unitality] \hfill\\ Writing $u := \otimes^0(*)$, and denoting $ab := a\otimes b = \otimes(a,b)$, the following must commute for every $a,b \in A$:

\begin{center}
{
\begin{tikzpicture}[node distance=1.7cm]

  \node (A1) {$ u u $} ;
  \node (A2) [right of=A1] {$ $};
  \node (B2) [below of=A2] {$u$};
  \node (B1) [below of=A1] {$ $};
  \node (B0) [left of=B1] {$u$};
  
  \draw[->, above right] (A1) to node {$ \lambda_u$} (B2);  
  \draw[->, above left] (B0) to node {$\rho_u$} (A1);
  \draw[->, below ] (B0) to node {$ 1_u $} (B2);
\end{tikzpicture}\hspace{.45cm}\begin{tikzpicture}[node distance=1.7cm]

  \node (A1) {$ab$} ;
  \node (A2) [right of=A1] {$ $};
  \node (B2) [below of=A2] {$u(a b)$};
  \node (B1) [below of=A1] {$ $};
  \node (B0) [left of=B1] {$(u  a) b$};
  
  \draw[->, above right] (B2) to node {$ \lambda_{a b}$} (A1);
  \draw[->, above left] (B0) to node {$ \lambda_a 1_b$} (A1);
  \draw[->, below ] (B0) to node {$ \alpha_{u,a,b} $} (B2);
\end{tikzpicture}

\medskip

\begin{tikzpicture}[node distance=1.7cm]

  \node (A1) {$ab$} ;
  \node (A2) [right of=A1] {$ $};
  \node (B2) [below of=A2] {$a(bu)$};
  \node (B1) [below of=A1] {$ $};
  \node (B0) [left of=B1] {$(ab)u$};
  
  \draw[->, above right] (A1) to node {$ 1_a \rho_b$} (B2);
  \draw[->, above left] (A1) to node {$ \rho_{ab}$} (B0);
  \draw[->, below] (B0) to node {$ \alpha_{a,b,u} $} (B2);
\end{tikzpicture} \hspace{.25cm}\begin{tikzpicture}[node distance=2.2cm]

  \node (A1) {$(au)b$} ;
  \node (A2) [right of=A1] {$a(ub)$};
  \node (B2) [below of=A2] {$ab$};
  \node (B1) [below of=A1] {$ab$};
  
  \draw[->, above] (A1) to node {$ \alpha_{a,u,b}$} (A2);  
  \draw[->, left] (B1) to node {$\rho_a 1_b $} (A1);
  \draw[->, right ] (A2) to node {$ 1_a \lambda_b $} (B2);
  
   \draw[-, line width=.08cm, right] (B1) to node {$  $} (B2); 
    \draw[-, white, line width=.03cm, below right] (B1) to node {$ $} (B2);

\end{tikzpicture}
}
\end{center}

}\end{description}

\end{definition}

\medskip

The authors of \cite{third} show that skew monoidal categories correspond precisely to maps out of $\ctn$ and into $\skn$, the \textit{skew nerve} of $\mathsf{Cat}$.\footnote{ It should be mentioned $\skn$ is referred to as the \textit{lax nerve} in \cite{third} which we have changed here so as to not suggest a relation to lax monoidal categories.}  This will be defined formally in section $5$.

\medskip

\begin{proposition}\cite{third}\label{garner}
Skew monoidal categories are in one to one correspondence with simplicial maps $\ctn \to\skn$.
\end{proposition}

\medskip

Lax monoidal categories are a second sort of weakened monoidal category characterized by the introduction of $n$-ary product operations.  The higher dimensional simplices of $\ctn$ capture these higher order products as well, as we will show shortly.

\medskip

\begin{definition}
A \textit{lax monoidal category} $(A, \otimes^n, \gamma, \iota)$ consists in:

\begin{enumerate} [\textbf{$\bullet$}]

\item A category $A$.
\item Functors $\otimes^n: A^n \to A$ for each $n\in \mathbb{N}$. ($\otimes^0: I \to A$) 
\item Natural transformations for each $n, k_1, ..., k_n \in \mathbb{N}$:
\[ \gamma_{n,k_1,...k_n} : \otimes^n \circ \left( \otimes^{k_1} \times ... \times \otimes^{k_n} \right) \tto \otimes^{k_1+...+k_n} \]
\item A natural transformation $\iota: 1_A \tto \otimes^1$.

\end{enumerate}

\medskip

\noindent  These natural transformations must additionally satisfy two axioms.

\begin{description}

\item[Associativity]\hfill\\  For each double sequence $n, k_1, ... ,k_n, m_{11},..., m_{1k_1},m_{21},...,m_{2k_2},...,m_{n1},...,m_{nk_n}$, the following square commutes:

\begin{center}\begin{tikzpicture}[node distance=3.5cm]

  \node (A1) {$ \otimes^n \circ (\otimes^{k_1} \times ... \times \otimes^{k_n} ) \circ ((\otimes^{m_{11}} \times... \times \otimes^{m_{1k_1}})\times...\times(\otimes^{m_{n1}} \times ... \times \otimes^{m_{nk_n}})) $} ;
  \node (A2) [right of=A1] {$ $};
  \node (B2) [below of=A2] {$\substack{\otimes^n \\ \circ (\otimes^{m_{11}+...+m_{1k_1}} \times ...\times \otimes^{m_{n1}+...+m_{nk_n}})}$};
  \node (B1) [below of=A1] {$ $};
  \node (B0) [left of=B1] {$\substack{\otimes^{k_1+...+k_n} \\ \circ ((\otimes^{m_{11}} \times... \times \otimes^{m_{1k_1}})\times...\times(\otimes^{m_{n1}} \times ... \times \otimes^{m_{nk_n}}))}$};
  \node (C1) [below of=B1] {$\otimes^{m_{11}+... +m_{1k_1}+...+m_{n1}+...+m_{nk_n}}$};
  
  \draw[->,  right] (A1) to node {$1 \circ (\gamma_{k_1,m_{11},...,m_{1k_1}} \times .... \times \gamma_{m_{n1},...,m_{nk_n}})$} (B2);  
  \draw[->,  below right] (B2) to node {$\gamma_{n,m_{11}+...+m_{1k_1},...,m_{n1}+...+m_{nk_n}}$} (C1);
  \draw[->, above left] (A1) to node {$\gamma_{n,k_1,...,k_n}\circ 1 $} (B0);
  \draw[->, below left] (B0) to node {$ \gamma_{k_1+...+k_n,m_{11},...,m_{1k_1},...,m_{n1},...,m_{nk_n} }$} (C1);
\end{tikzpicture}\end{center}

\item[Unitality]\hfill\\  The following two triangles commute:

\begin{center}\begin{tikzpicture}[node distance=2.5cm]

  \node (A1) {$ \hspace{1.85cm} 1_A \circ \otimes^n = \otimes^n = \otimes^n \circ (1_A \times... \times 1_A)$} ;
  \node (A2) [right of=A1] {$ $};
  \node (B2) [below of=A2] {$\otimes^n \circ (\otimes^1 \times ... \times \otimes^1) $};
  \node (B1) [below of=A1] {$ $};
  \node (B0) [left of=B1] {$\otimes^1 \circ \otimes^n $};
  \node (C1) [below of=B1] {$\otimes^n $};
  
  \draw[->,  above right] (A1) to node {$ 1 \circ \iota \times ... \times \iota $} (B2);  
  \draw[->,  below right] (B2) to node {$\gamma_{n,1,...,1}$} (C1);
  \draw[->, above left] (A1) to node {$\iota \circ 1 $} (B0);
  \draw[->, below left] (B0) to node {$ \gamma_{1,n} $} (C1);
  \draw[->, below left] (A1) to node {$ 1 $} (C1);
\end{tikzpicture}\end{center}

\end{description}

\end{definition}

\medskip

In particular, the data of a lax monoidal category includes a binary functor $\otimes^2:A\times A \to A$.  Though it is not required to be associative or even associative up to natural isomorphism, we are given a pair of natural transformations:

\[ \gamma_{2,2,1} \bullet (1 \circ (1 \times \iota)) : \otimes^2 \circ (\otimes^2 \times 1) \tto \otimes^3\]
\[ \gamma_{2,1,2} \bullet (1 \circ (\iota \times 1)) : \otimes^2 \circ (1 \times \otimes^2) \tto \otimes^3\]

\noindent where we use `$\bullet$' to denote vertical composition.  Writing $\otimes^n(a_1,...,a_n)= (a_1 \otimes ... \otimes a_n)$, the above natural transformations give us maps:

\[ ((a \otimes b) \otimes c)) \rightarrow (a \otimes b \otimes c) \leftarrow (a \otimes (b \otimes c)) \]

\noindent  In this way, the $n$-ary product operations mediate between composites of lower arity operations, and hence say something about a weakened form of associativity.

\medskip

Whereas skew monoidal categories can be understood as maps $\ctn \to \skn$, we will see that lax monoidal categories can be understood as maps $\ctn \to \hcn$, a nerve containing $\skn$ in a way we will make precise in section $5$.  We turn now to define this nerve.

\medskip

\begin{definition}
The \textit{homotopy coherent nerve} $N_{\Delta}(U)$ of a simplicially enriched category $U$ is the simplicial set defined on objects by $N_\Delta(U)_n := \mathsf{sSetCat}(S[n],U)$.  The enriched category $S[n]$ has as objects the elements $0,1,..., n \in \n$.  For each pair of objects $p \le q$, the simplical mapping object $S[n](p,q)$ is the categorical nerve of the poset whose objects are subsets of $\lbrace p , p+1 ,...  , q\rbrace$ containing both $p$ and $q$ and ordered by inclusion.  (For $q < p$, $S[n](p,q) = \emptyset$.)  Here are some examples of simplices in these simplicial mapping objects:

\begin{enumerate}[$\bullet$]

\item $\left( \lbrace 0, 3 \rbrace \right) \in S[3](0,3)_0$
\item $\left(\lbrace 0,3 \rbrace \subset [3] \right) \in S[3](0,3)_1$
\item $\left(\lbrace 0,3 \rbrace \subset \lbrace 0,1,3 \rbrace \subset [3]\right) \in S[3](0,3)_2$
\item $\left(\lbrace 0,3 \rbrace \subset \lbrace 0,1,3 \rbrace= \lbrace 0,1,3 \rbrace \right) \in S[3](0,3)_2$
\item $\left(\lbrace 1,4,5 \rbrace \subset \lbrace 1,2,4,5 \rbrace \right) \in S[9](1,5)_1$

\end{enumerate}

\noindent We will be careful to reserve the symbol `$\subset$' to mean a \textit{proper} inclusion, and will use `$\subseteq$' otherwise. 

Composition maps $S[n](p,r) \times S[n](r,q) \to S[n](p,q)$ are given by the union of subsets in dimension $0$, and unions of inclusions of subsets in higher dimensions.  The identity elements are therefore $\lbrace p \rbrace \in S[n](p,p)_0$.  $N_\Delta(U)_n$ then consists of all simplicially enriched functors out of $S[n]$ and into $U$.  

\medskip

The coface and codegeneracy maps of $\Delta$ extend to enriched coface functors $  \delta_i:S[n-1] \to S[n]$ and codegeneracy functors $\sigma_i: S[n+1] \to S[n]$.  On objects, these functors match their counterparts in $\Delta$,  while on mapping objects, $\delta_i:S[n-1](p,q) \to S[n](\delta_i p, \delta_i q)$ sends a subset $\ccc$ in dimension $0$ to its direct image $\delta_i \ccc = \lbrace \delta_i c_0, ... ,\delta_i c_m \rbrace$, and sends inclusions to their direct images in higher dimensions. Codegeneracy functors are defined similarly.  Finally, the face and degeneracy maps of $ N_\Delta(U)$ are defined via precomposition with these enriched functors: $d_i : N_\Delta(U)_n \to N_\Delta(U)_{n-1}$ is given by $  d_i(L) = L \circ \delta_i $ for $L \in N_{\Delta}(U)_n$, and similarly for $s_i$. 
\end{definition}

\medskip

In what follows, we will focus specifically on the homotopy coherent nerve of $\mathsf{Cat}$, viewed as a simplicially enriched category in the following way. There is a standard nerve functor $N_2: \mathsf{2Cat} \to \mathsf{sSet}$ which preserves products, thus taking monoids in $\mathsf{2Cat}$ to monoids in $\mathsf{sSet}$.  Monoids in $\mathsf{sSet}$ can in turn be viewed as one object simplicially enriched categories with the monoid as the lone simplicial mapping object.  In summation, we have a nerve $N_2: \mathsf{Mon}_{\mathsf{2Cat}}\to \mathsf{sSetCat}$ which takes a monoid in $\mathsf{2Cat}$ to a one object simplicially enriched category. We would like to apply this nerve to $\mathsf{Cat}$, an object of $\mathsf{2Cat}$.  Unfortunately, $\mathsf{Cat}$ is not a monoid, strictly speaking.  With $\times$ as a binary operation and the terminal category $I$ as unit, $\mathsf{Cat}$ is itself a monoidal category, a not-quite-monoid of $\mathsf{2Cat}$: the operation $\times$ is associative and unital merely up to natural isomorphism.  However, according to the coherence theorem \cite{second}, we know that $\mathsf{Cat}$ is monoidally equivalent to another $2$-category which is an actual monoid.  By $N_2(\mathsf{Cat})$ we mean $N_2$ applied to this equivalent category.  In practice, we consider the following definition.

\medskip

\begin{definition}

The simplicially enriched category $N_2(\mathsf{Cat})$ is defined by the following data:

\begin{enumerate}[\textbf{$\bullet$}]

\item A unique object, $*$.

Its mapping object $N_2(\mathsf{Cat})(*,*)$ is characterized by:

\item $0$-simplices are Categories $B$.
\item $1$-simplices are functors $T: B_1 \to B_0$.
\item $2$-simplices are natural transformations $\eta: T_{12}\circ T_{01} \tto T_{02}$.
\item $3$-simplices are commutative diagrams of natural transformations:

\begin{center}\begin{tikzpicture}[node distance=2.5cm]
 \node (A1) {$ T_{23} \circ T_{12} \circ T_{01}$} ;
  \node (A2) [right of=A1] {$ $};
  \node (B2) [below of=A2] {$T_{13}\circ T_{01}$};
  \node (B1) [below of=A1] {$ $};
  \node (B0) [left of=B1] {$T_{23} \circ T_{02}$};
  \node (C1) [below of=B1] {$T_{03}$};
  
  \draw[->,  right] (A1) to node {$ 1 \circ \eta_{123}$} (B2);  
  \draw[->,  below right] (B2) to node {$\eta_{013}$} (C1);
  \draw[->, above left] (A1) to node {$ \eta_{012} \circ 1 $} (B0);
  \draw[->, below left] (B0) to node {$ \eta_{023}$} (C1);
\end{tikzpicture}\end{center}

\item Higher-dimensional simplices are given by $3$-coskeletality.  

\end{enumerate}

Composition in $N_2(\mathsf{Cat})$ is given by $\times$ while the identity map for composition consists in a $0$-simplex of $N_2(\mathsf{Cat})(*,*)$: $I$.  Face maps in $N_2(\mathsf{Cat})(*,*)$ are given as suggested by the notation above, \\ e.g $d_1 (\eta:T_{12} \circ T_{01} \to T_{02}) = T_{02}$, and degeneracy maps are given by inserting identity maps in the expected ways.  We assume that we have the following equalities: 

\[ (A\times B) \times C = A \times (B \times C) \]
\[ A \times I = A = I \times A \]

\end{definition}

It is worth noting that $3$-coskeletality of $N_2(\mathsf{Cat})(*,*)$ shows us that every simplex is determined by its $3$-faces, but also we see that the $3$-simplices are determined by their $3$-boundaries: there is at most one $3$-simplex with a given $3$-boundary.  Thus, every simplex of $N_2(\mathsf{Cat})(*,*)$ is determined by its $2$-faces.

\medskip

We can give an informal description of simplices of $\hcn$.  It has a single $0$-simplex: the unique map $S[0] \to N_2(\mathsf{Cat})$.  Its $1$-simplices can be thought of as categories $A$, while its $2$-simplices consist in functors $T:A_{01} \times A_{12} \to A_{02}$.  A $3$-simplex $L:S[3] \to  N_2(\mathsf{Cat})$ is a diagram consisting of five such functors, four of which are the $3$-simplex's faces.  We can see this diagram most clearly by thinking of the image of $L$ on the simplicial mapping object $S[3](0,3)$.  $S[3](0,3)$ has four $0$-simplices (subsets of $[3]$ containing $\lbrace 0,3 \rbrace$), five $1$-simplices (inclusions of those subsets), and two $2$-simplices (double inclusions).  The image therefore consists in four categories, five functors, and two natural transformations.  This data fits into the following diagram:

\begin{center} \begin{tikzpicture}[node distance=2cm]
  \node (A1) {$\lbrace 0,1,2,3 \rbrace$} ;
  \node (A2) [right of=A1] {$ $};
  \node (A3) [right of=A2] {$ \lbrace 0,1,3\}$};
  \node (B2) [below of=A2] {$ $};
  \node (B1) [below of=A1] {$ $};
  \node (B3) [below of=A3]{$ $};
  \node (C1) [below of=B1] {$\{0,2,3\}$};
  \node (C2) [below of=B2] {$ $};
  \node (C3) [below of=B3] {$\{0,3\}$};
  
  \draw[->,  above] (A1) to node {$ $} (A3);  
  \draw[->,  left] (A1) to node {$ $} (C1);
  \draw[->, below] (C1) to node {$ $} (C3);
  \draw[->, right] (A3) to node {$ $} (C3);
  \draw[->, above left] (A1) to node {$ $} (C3);

\end{tikzpicture}\hspace{.6cm}\begin{tikzpicture}[node distance=2cm]
  \node (A1) {$ A_{01}\times A_{12} \times A_{23}$} ;
  \node (A2) [right of=A1] {$ $};
  \node (A3) [right of=A2] {$ A_{01}\times A_{13}$};
  \node (B2) [below of=A2] {$ $};
  \node (B1) [below of=A1] {$ $};
  \node (B3) [below of=A3]{$ $};
  \node (C1) [below of=B1] {$A_{02}\times A_{23}$};
  \node (C2) [below of=B2] {$ $};
  \node (C3) [below of=B3] {$A_{03}$};
  
  \draw[->,  above] (A1) to node {$1_{A_{01}} \times T_{123}$} (A3);  
  \draw[->,  left] (A1) to node {$T_{012} \times 1_{A_{23}}$} (C1);
  \draw[->, below] (C1) to node {$T_{023}$} (C3);
  \draw[->, right] (A3) to node {$T_{013}$} (C3);
  \draw[->,  left] (A1) to node {$ $} (C3);
  \draw[->, below left] (A3) to node {$ $} (B2);
 \draw[->, line width=.08cm, above left] (A3) to node {$   \eta_{\lbrace 0,1,3 \rbrace}$} (B2); 
 \draw[-, white, line width=.03cm, below right] (A3) to node {$ $} (B2);   
 \draw[->, line width=.08cm, below right] (C1) to node {$  \eta_{\lbrace 0,2,3 \rbrace}$} (B2); 
 \draw[-, white, line width=.03cm, below right] (C1) to node {$ $} (B2);   

\end{tikzpicture}\end{center}
\[ \hspace{.8cm} S[3](0,3) \hspace{5.8cm} L(S[3](0,3)) \]

\medskip\medskip

Higher dimensional simplices $L:S[n] \to N_2(\mathsf{Cat})$ are again helpfully summarized via their image of the simplicial mapping object $S[n](0,n)$.  Such a simplex will provide: a category for each subset of $\n$ containing $\lbrace 0,n \rbrace$; a functor between such categories whenever the subset indexing the first contains the subset indexing the second; a natural transformation between a composite of such functors and a third such functor for each double containment; a commutative diagram of such natural transformations for every triple and higher containment of subsets.  If $\ccc \in S[n](0,n)_0$, because $\ccc = \cup_{c \in \ccc^-}\lbrace c,sc \rbrace$ and $L$ must send unions to products, we have then that $L(\ccc) = \prod_{c \in \ccc^-} L(\lbrace c,sc \rbrace)$, where the right hand side comes from $L$'s action on the simplicial mapping objects $S[n](c,sc)$.  Functors and natural transformations may be given as products in this way as well.  Notably, the images of the simplices of $S[n](0,n)$ of the form $\lbrace 0,n \rbrace \subseteq \n$ and $\lbrace 0,n \rbrace \subseteq \ccc \subseteq \n$ are of particular importance: the image of every other $1$ and $2$-simplex of any simplicial mapping object $S[n]$ appears as the data of a proper face of $L$.  For example, if $L$ is a $3$-simplex, we can see from the above that the functor $L(\lbrace 0,3 \rbrace \subset \lbrace 0,2,3 \rbrace)$ is actually the functor associated with the $2$-simplex $d_1(L)$.  This phenomena continues into the higher dimensions.   In what follows, we will study the structure of $\hcn$ in much greater detail and rigour.

 We are now ready to state and prove the three main results.  In section $4$ of this paper we will, in the spirit of Proposition \ref{garner}, prove:

\begin{proposition}\label{laxclass}
There is an assignment of a map $\alpha : \ctn \to N_\Delta (N_2\mathsf{Cat})$ to each lax monoidal $(A,\otimes^n, \iota, \gamma)$ such that $(A,\otimes^n, \iota,\gamma)$ can be recovered completely from $\alpha$.
\end{proposition}

\medskip

In section $5$ we will extend and unify the previous results concerning skew monoidal categories by showing that:

\begin{proposition}\label{skewclass}
There is a monomorphism $\beta:\skn \hookrightarrow \hcn$.
\end{proposition}

\medskip

\noindent Combined with Proposition \ref{garner}, we will then have that skew monoidal categories correspond to maps $\ctn \to \hcn$ as well.  

In section $6$ we will define $\Sigma$-monoidal categories for a countable signature $\Sigma$, and prove:

\begin{proposition}\label{sigmaclass}
There is an assignment of a map $\sigma : \ctn \to \hcn$ to each $\Sigma$-monoidal category $(A, \Sigma, \gamma)$ such that $(A, \Sigma, \gamma)$ can be recovered completely from $\sigma$. \end{proposition}

\medskip

Taken together, our three results along with the crucial proposition \ref{backward} characterizing maps $\ctn \to \hcn$ show that maps $\ctn \to \hcn$ are a natural setting to understand a great many monoidal-type categories, including ones not yet defined in the literature.  As a first and significant step, we will examine arbitrary maps $\ctn \to \hcn$ in detail and identify an alternative way of defining them in terms of a simple list of data subject to some few coherence requirements.

\section{Defining maps into $\hcn$}

Defining maps into $\hcn$ is not as daunting as it might first appear.  In this section, we will explore these maps and develop a succinct way of producing them.   In what follows, let $X$ be an arbitrary simplicial set.

\begin{proposition}\label{forward}

A map $\phi: \ds \to \hcn$ is determined by the values:

\begin{enumerate}
\item $\phi_1(x)(\lbrace 0,1 \rbrace)$ for each non-degenerate $x \in \ds_1$.

\item $\phi_n(x)(\lbrace 0,n \rbrace \subset \n)$ for each non-degenerate $x \in \ds_n$, $n \ge 2$.

\item $\phi_n(x)(\lbrace 0,n \rbrace \subset \ccc \subset \n)$  for each non-degenerate $x \in \ds_n$, $n \ge 3$, and non-degenerate $\left(\lbrace 0,n \rbrace \subset \ccc \subset \n\right) \in S[n](0,n)_2$.
\end{enumerate}

Moreover, these values can be explicitly described:

\begin{enumerate}
\item  $\phi_1(x)(\lbrace 0,1 \rbrace) = A^x$, a category.
\item $\phi_n(x)(\lbrace 0,n \rbrace \subset \n) = T^x: \prod \limits_{i \in \n^-} A^{x_{i,i+1}} \to A^{x_{0,n}}$, a functor.
\item $\phi_n(x)(\lbrace 0,n \rbrace \subset \ccc \subset \n) = \eta^x_{\ccc}: (T^{x_\ccc} \circ \prod\limits_{c \in \ccc^-}T^{x_{[c,sc]}}) \tto T^x$, a natural transformation. \footnote{  Here it may very well be that any of $x_{i,i+1}$, $x_{\ccc}$, or $x_{[c,sc]}$ is degenerate, or even that $ [c,sc]=\{c,sc\}$ so that $x_{c,sc}$ is a $1$-simplex.  By $A^y$ and $T^y$ we mean always $\phi_n(y)(\lbrace 0,1 \rbrace)$ and $\phi_n(y)(\{0,n\}\subseteq [n])$ respectively, even when $y$ is degenerate or $n=1$.}

\end{enumerate}

\end{proposition}

\medskip

\begin{proof}

Suppose that $\phi:\ds \to \hcn$ and that we know the values specified in items $(1), (2)$, and $(3)$ above.  We will show that all other values of $\phi$ can be determined from these.  

\medskip

For $x \in \ds_n$,  $\phi_n(x):S[n] \to N_2(\mathsf{Cat})$ is trivial on objects because $N_2(\mathsf{Cat})$ has a single object $*$.  We consider then its action on the mapping spaces $S[n](p,q)$ of $S[n]$. For brevity, we call a $k$-simplex of a mapping space of $S[n]$ a $k$-simplex of $S[n]$.  Note that if $x\in \ds$ is degenerate, because $\phi$ commutes with degeneracy maps, its value on $x$ is determined by its value on the unique non-degenerate simplex mapping to $x$, hence it suffices to consider non-degenerate simplices.  In what follows, we suppose $x$ is always non-degenerate.

As $\phi_n(x)$ is an enriched functor, it must respect composition in $S[n]$.  For $\ccc$ a $0$-simplex of $S[n]$, we have: 
\[\phi_n(x)(\ccc)= \phi_n(x)\left(\bigcup_{c\in \ccc^-} \lbrace c ,sc \rbrace\right) = \prod \limits_{c \in \ccc^-} \phi_n(x)(\lbrace c,sc \rbrace )\]

This shows that $\phi_n(x)$ is determined on $0$-simplices of $S[n]$ by its values on subsets of the form $\lbrace p,q \rbrace$ with $0 \le p < q \le n$. For each such $\lbrace p,q \rbrace$ there is the coface map $\delta_{p,q}:S[1] \to S\n$ sending $0$ to $p$ and $1$ to $q$.  This gives:

\[ \phi_n(x)(\lbrace p,q \rbrace) = \phi_n(x)(\delta_{p,q}\lbrace 0,1 \rbrace) = \phi_n(x) \circ \delta_{p,q} (\lbrace 0,1 \rbrace) = \phi_1(x_{p,q})(\lbrace 0,1\rbrace)\]

The last equality follows from $\phi$ commuting with face maps. If $x_{p,q}$ is degenerate, then because $\phi$ commutes with degeneracy maps, we must have $\phi_1(x_{p,q})(\lbrace 0,1 \rbrace) = I$ while if $x_{p,q}$ is non-degenerate, its value is specified by item $(1)$.  Hence $\phi_n(x)$ is specified on all $0$-simplices by the data of item $(1)$.

\medskip

Again because $\phi_n(x)$ respects composition, $\phi_n(x)$ is determined on arbitrary $1$-simplices $(\ccc_0 \subseteq \ccc_1)$ by $1$ simplices of the form $(\lbrace p,q \rbrace \subseteq \ccc)$.  Because $\phi$ commutes with face maps, we have:

\[ \phi_n(x)(\lbrace p,q \rbrace \subseteq \ccc) = \phi_n(x)(\delta_{\ccc}(\lbrace 0,m \rbrace \subseteq [m])) = \phi_m(x_{\ccc})(\lbrace 0,m\rbrace \subseteq [m])\]

\noindent Such data is specified in item $(2)$ in the case $x_{\ccc}$ is non-degenerate.  If however $x_{\ccc}$ is the degeneracy of some non-degenerate simplex $y$ of dimension $l$, $x_{\ccc}=sy$, we have:

\[\phi_m(x_{\ccc})(\lbrace 0,m\rbrace \subseteq [m]) = \phi_m(sy)(\lbrace 0,m\rbrace \subseteq [m]) = \phi_l(y)(\lbrace \sigma 0, \sigma m\rbrace \subseteq \sigma[m])= \phi_l(y)(\lbrace 0,l\rbrace \subseteq [l]) \]

\noindent Thus this is determined again by the data of item (2).\footnote{In the case that $\sigma 0 = 0$, $\sigma m = 1$, then $y$ is a $1$-simplex and $\phi_1(y)(\{0,1\} \subseteq [1]) = 1_{\phi_1(y)(\{0,1\})}$}.

\medskip

Similarly, the value of $\phi_n(x)$ on $2$-simplices is determined by its values on simplices $(\lbrace 0,n \rbrace \subseteq \ccc \subseteq \n)$ as a result of both respecting composition and commuting with face maps.  It is therefore determined by item $(3)$, noting that we may need to commute with degeneracy maps as above. Finally, because $\phi_n(x)(\ccc_0 \subseteq \ccc_1 \subseteq ... \subseteq \ccc_k)$ must be a $k$-simplex of $N_2(\mathsf{Cat})$, it is hence determined by its $2$-faces.  Because $\phi_n(x)$ commutes with face and degeneracy maps, these $2$-faces are determined by the value of $\phi_n(x)$ on the $2$-faces of ($\ccc_0 \subseteq \ccc_1 \subseteq ... \subseteq \ccc_k$) and hence by the data of item (3) as above.

\medskip

As for the explicit descriptions, $\phi_1(1)(\lbrace 0,1 \rbrace)$ is a category,  $\phi_n(x)(\lbrace 0,n \rbrace \subset \n)$, a functor, and $\phi_n(x)(\lbrace 0,n \rbrace \subset \ccc \subset \n)$ a transformation from a composite of functors, all in light of the definition of $N_2(\mathsf{Cat})$.  We know that the codomain of the functor $T^x$ is :  

\[d_0 T^x = d_0 \phi_n(x)(\lbrace 0,n \rbrace \subset \n) = \phi_n(x)(d_0 (\lbrace 0,n \rbrace \subset [n])) = \phi_n(x)(\lbrace 0,n \rbrace) = \phi_1(x_{0,n})(\lbrace 0,1 \rbrace) = A^{x_{0,n}}\]  

\medskip

\noindent Its domain is:

\[d_1 T^x = d_1\phi_n(x)(\lbrace 0,n \rbrace \subset \n) = \phi_n(x)(d_1 (\lbrace 0,n \rbrace \subset [n])) = \phi_n(x)([n]) = \prod \limits_{i \in [n]^-} \phi_1(x_{i,i+1})(\lbrace 0,1 \rbrace) = \prod \limits_{i \in \n^-}A^{x_{i,i+1}}\]

\noindent  We can see the codomain of $\eta^x_{\ccc}$  by considering:

\[d_1 \eta^x_{\ccc} =  d_1 \phi_n(x)( \lbrace 0,n \rbrace \subset \ccc \subset [n]) =  \phi_n(x)(d_1 (\lbrace 0,n \rbrace \subset \ccc \subset [n])) = \phi_n(x)( \lbrace 0,n \rbrace \subset [n]) = T^x \]

\medskip

\noindent  And can see its domain by considering both of:

\[d_0 \eta^x_{\ccc} = \phi_n(x)(d_0 (\lbrace 0,n \rbrace \subset \ccc \subset [n])) = \phi_n(x)( (\lbrace 0,n \rbrace \subset \ccc)) = T^{x_{\ccc}} \]
\[d_2 \eta^x_{\ccc} = \phi_n(x)(d_2 (\lbrace 0,n \rbrace \subset \ccc \subset [n])) = \phi_n(x)( \ccc \subset [n]) = \prod \limits_{c \in \ccc^-} \phi_n(x)( \lbrace c,sc \rbrace \subseteq [c,sc]) = \prod \limits_{c \in \ccc^-} T^{x_{[c,sc]}}\]

\end{proof}

\medskip

There is also a converse to the above proposition which tells us exactly the conditions a collection of categories, functors, and natural transformations needs to satisfy in order to extend to a map $X \to \hcn$.  We state it properly (proposition \ref{backward}) and prove it below, but will first need some key observations about the homotopy coherent nerve construction.  Let $U$ be a simplicially enriched category.  We know that $n$-simplices of $N_{\Delta}(U)$ are given by simplicially enriched functors $S[n] \to U$.  The mapping simplicial sets of $S[n]$ are freely generated via union by simplices of the form $\lbrace p, q \rbrace$, $\lbrace p,q \rbrace \subseteq \ccc_1$, $\lbrace p,q \rbrace \subseteq \ccc_1 \subseteq \ccc_2$, and so on, as we have seen in the preceding proof.  Therefore a map $L:S[n] \to U$ is determined by where it sends the objects $0,1,...,n$, and where it sends the generating simplices of the mapping sets just listed.  If we are in the situation that the simplicial mapping objects of $U$ are uniformly $r$-coskeletal for some $r \ge 0$, then $L:S[n] \to U$ is determined by where it sends the objects and generating simplices of dimension $\le r$, that is, generating simplices up to those of the form $\lbrace p,q \rbrace \subseteq \ccc_1  \subseteq ... \subseteq \ccc_r$.  

\medskip

The simplicially enriched category $N_2(\mathsf{Cat})$ has a $3$-coskeletal simplicial mapping object.  Indeed, it is `nearly' $2$-coskeletal, as every $3$-boundary is the boundary of at most one $3$-simplex.  Writing it all out explicitly, we get the following lemma. 

\medskip

\begin{lemma}\label{Llemma}

To define a simplicially enriched functor $L:S[n] \to N_2(\mathsf{Cat})$, it suffices to define

\begin{enumerate}

\item A category $L(\lbrace p, q \rbrace)$ for all $0 \le p < q \le n$
\item A functor 
\[ L(\lbrace p,q \rbrace \subseteq \ccc):  \prod \limits_{c\in \ccc^-} L(\lbrace c,sc \rbrace) \to L(\lbrace p,q \rbrace)\] for each $\lbrace p,q \rbrace \subseteq \ccc $ such that $ L(\{p,q\} = \{p,q\})= 1_{L(\{p,q\})}$

\item A transformation \[L(\lbrace p,q \rbrace \subseteq \ccc_1 \subseteq \ccc_2): L(\lbrace p,q \rbrace \subseteq \ccc_1) \circ \prod \limits_{c \in \ccc_1^-} L(\lbrace c,sc \rbrace \subseteq \ccc_2 \cap [c,sc]) \tto L(\lbrace p,q \rbrace \subseteq \ccc_2)\]for each $\lbrace p,q \rbrace \subseteq \ccc_1 \subseteq \ccc_2$  such that $L(\lbrace p,q \rbrace \subseteq \ccc_1 \subseteq \ccc_2)= 1_{L(\{p,q\}\subseteq \ccc_2)}$ if $\ccc_1 = \{p,q\}$ or $\ccc_1 = \ccc_2$.

\end{enumerate}

\medskip \vspace{.3cm}

And such that, for every non-degenerate $(\{p,q\}\subset \ccc_1 \subset \ccc_2 \subset \ccc_3) \in S[n](p,q)_3$:

\begin{align*}\label{key2}\tag*{$(\star)$} &L(\{p,q\} \subset \ccc_2 \subset \ccc_3) \bullet (L( \{p,q\} \subset \ccc_1 \subset \ccc_2) \circ 1) = \\& L(\{p,q\} \subset \ccc_1 \subset \ccc_3) \bullet \left( 1 \circ \prod \limits_{c \in \ccc_1^-} L(\{c,sc\} \subseteq \ccc_2\cap [c,sc] \subseteq \ccc_3 \cap [c,sc])\right)
\end{align*}

\end{lemma}

Here the specification of domains and codomains of $ L(\lbrace p,q \rbrace \subseteq \ccc)$ and $L(\lbrace p,q \rbrace \subseteq \ccc_1 \subseteq \ccc_2)$ is precisely what is required for $L$ restricted to simplicial mapping objects of $S[n]$ to commute with face maps.  The qualification that, e.g. $L(\lbrace p,q \rbrace = \lbrace p,q \rbrace) = 1_{L(\lbrace p,q \rbrace)}$ is what is required for $L$ restricted to simplicial mapping objects to commute with degeneracy maps.  Finally, had the mapping object of $N_2(\mathsf{Cat})$ truly been $2$-coskeletal, we could have done away with equation \ref{key2}.  Be that as it may, the $3$-boundary consisting of the four natural transformations appearing in equation \ref{key2} is the boundary of a $3$-simplex precisely when those transformations commute, i.e, when \ref{key2} is satisfied.  Finally, if any of the `$\subset$'s of $(\{p,q\}\subset \ccc_1 \subset \ccc_2 \subset \ccc_3) $ had been `$=$'s, this commutativity would have been given automatically explaining why we needn't consider such cases.

\medskip

Proposition \ref{backward} -- the converse to proposition \ref{forward} -- takes a bit of work to state correctly.  Towards these ends, we will need the following three technical lemmas.

\begin{lemma}\label{degens}

Let $\ds$ be a simplicial set, $\sigma : [n]\to [n-1]$ an arbitrary codegeneracy map with corresponding degeneracy $s$.  Let $\ccc \subseteq [n]$ with $|\ccc|=m+1$.  Then for all $y\in \ds_{n-1}$, $(sy)_{\ccc}=s'(y_{\sigma \ccc})$ where $s'$ is a degeneracy map or an identity, and is the latter if and only if $|\ccc| = |\sigma \ccc| = m+1$.

\end{lemma}

\medskip

\begin{proof}

 Let $\delta_{\ccc}: [m] \to \n$ be the composition of coface maps sending $i \mapsto c_i$.  We have that $\sigma \delta_{\ccc}$ can be rewritten as a surjection followed by an injection.  We know that $|\sigma \delta_{\ccc} [m]| = |\sigma \ccc| = m+1 $ or $m$, which tells us that the surjection in the rewrite either has target $[m]$ or $[m-1]$, hence is an identity or a codegeneracy map which we denote in either case as $\sigma'$.  So we have  $\sigma \delta_{\ccc} = \delta_{\sigma \ccc} \sigma'$.  In terms of faces and degeneracies we then have $d_{\ccc} s = s' d_{\sigma \ccc}$ and so have:
 
 \[ (sy)_{\ccc} = d_{\ccc} s y = s' d_{\sigma \ccc} y = s'(y_{\sigma \ccc}) \]

\end{proof}

\medskip

\begin{lemma}\label{part1}
Let $A^x$ be a category for each $x \in \ds_1$ such that $A^{sy} = I$ for the degeneracy of a $0$-simplex $y$.  For every $n \ge 1$ and $x \in \ds_n$, define the symbols:
\[ \dom(x):= \prod \limits_{i \in \n^-} A^{x_{i,i+1}}  \hspace{.4cm} , \hspace{.4cm} \cod(x):= A^{x_{0,n}}\]
\noindent Then we have:

\begin{enumerate}[(i)]
\item $\dom(x)=A^x= \cod(x)$ if $x \in \ds_1$.
\item $\dom(sx)=\dom(x)$ and $\cod(sx)=\cod(x)$ for any degeneracy map $s$.

And for every $\lbrace 0,n \rbrace \subseteq \ccc \subseteq \n$,

\item $\dom(x_{\ccc})=\prod \limits_{c \in \ccc^-} A^{x_{c,sc}}$ and $\cod(x_{\ccc})=A^{x_{c_0,c_m}}$.
\end{enumerate}

\end{lemma}

\medskip

\begin{proof}

We consider each item in turn.

\begin{enumerate}[(i)]

\item If $x \in \ds_1$, then as $[1]^- = \lbrace 0 \rbrace$, $\dom(x) = A^{x_{0,1}} = A^x = \cod(x)$.

\item If $sx$ is the degeneracy of a simplex $x \in \ds_n$, then by lemma \ref{degens}, $(sx)_{i,i+1}=s'(x_{\sigma i, \sigma i+1})$ where $s'$ is the degeneracy map $\ds_0 \to \ds_1$ if $\sigma i = \sigma( i+1)$, i.e if $\sigma = \sigma_i$, and is otherwise the identity.  Suppose then that $\sigma= \sigma_k$ for some $k \in \n$.  Then we have:

\begin{align*}
\dom(sx)= & A^{s'(x_{\sigma 0, \sigma 1})} \times A^{s'(x_{\sigma 1, \sigma 2})} \times ... \times A^{s'(x_{\sigma k, \sigma k+1})} \times ... \times A^{s'(x_{\sigma n-1, \sigma n})} \times A^{s'(x_{\sigma n, \sigma n+1})} \\
 = & A^{x_{ 0, 1}} \times A^{x_{1, 2}} \times ... \times A^{s'(x_{ k, k})} \times ... \times A^{x_{ n-2,  n-1}} \times A^{x_{n -1, n}}\\
= & A^{x_{ 0, 1}} \times A^{x_{1, 2}} \times ... \times I \times ... \times A^{x_{ n-2,  n-1}} \times A^{x_{n -1, n}} \\
=& \dom(x)
\end{align*}

\medskip

\noindent Similarly, $\cod(sx)=(sx)_{0,n+1} = s'(x_{\sigma 0,\sigma n+1}) = x_{0,n} =\cod(x)$ because $s'$ is always the identity (as $n \ge 1$).

\item Given $x \in \ds_n$ and $\lbrace 0,n \rbrace \subseteq \ccc \subseteq \n$, we have $x_{c_i,c_{i+1}} = (x_{\ccc})_{i,i+1}$ simply by composing the face maps $\delta_{\ccc}$ and $\delta_{i,i+1}$.  The result follows.

\end{enumerate}

\end{proof}

\medskip

\begin{lemma}\label{part2}
Let $A^x$ be a category for each $x \in \ds_1$ such that $A^{sy} = I$ for the degeneracy of a $0$-simplex $y$.  For each $n \ge 1$ and $x \in \ds_n$, let $T^x: \dom(x) \to \cod(x)$ be a functor such that $T^x = 1_{A^x}$ if $x \in \ds_1$ and $T^{sx}=T^x$.  For every $n \ge 1$, $x \in \ds_n$, and $\lbrace 0,n \rbrace \subseteq \ccc \subseteq \n$, define the symbols:

\begin{center} $\dom(x,\ccc) := T^{x_{\ccc}} \circ \prod \limits_{c \in \ccc^-} T^{x_{[c,sc]}} \hspace{.4cm} , \hspace{.4cm} \cod(x,\ccc) := T^x $  \footnote{The composite of $\dom(x,\ccc)$ is well defined by item (iii) of lemma \ref{part1}}\end{center}

\noindent Then we have:

\begin{enumerate}[(i)]
\item $\dom(x,\ccc)=T^x=\cod(x,\ccc)$ if $\ccc = \lbrace 0,n \rbrace$ or $\ccc = \n$.
\item $\dom(sx,\ccc) = \dom(x,\sigma \ccc)$ and $\cod(sx,\ccc)=\cod(x,\sigma \ccc)$ for any degeneracy map $s$.

And for every $\lbrace 0,n \rbrace \subseteq \aaa \subseteq \bbb \subseteq \n$,

\item $\dom(x_{\bbb}, \delta_{\bbb}^{-1}\aaa) = T^{x_{\aaa}} \circ \prod \limits_{a \in \aaa^-} T^{x_{\bbb \cap [a,sa]}}$ and $\cod(x_{\bbb}, \delta_{\bbb}^{-1}\aaa) = T^{x_{\bbb}}$.

\end{enumerate}

\end{lemma}

\medskip

\begin{proof} We consider each item in turn.

\begin{enumerate}[(i)]

\item If $\lbrace 0,n \rbrace \subseteq \ccc \subseteq \n$ and $\ccc = \{0,n\}$, then $x_{\ccc} = x_{0,n}$ is a $1$-simplex, hence $T^{x_{\ccc}} = 1$ and $\dom(x,\ccc) = 1 \circ T^{x_{[0,n]}} = T^x$.  If $\ccc = \n$, then $x_{\ccc} = x$, and $T^{x_{[c,sc]}} = T^{x_{i,i+1}}$ for some $i$, and is hence an identity.  Thus $\dom(x,\ccc) = T^x \circ \prod 1 = T^x$.

\item If $sx$ is the degeneracy of a simplex $x \in \ds_n$ and  $\lbrace 0,n \rbrace \subseteq \ccc \subseteq \n$, then by lemma \ref{degens}, $(sx)_{\ccc} = s'(x_{\sigma \ccc})$ and $(sx)_{[c,sc]} = s'(x_{[\sigma c, \sigma sc]})$ so that $T^{(sx)_{\ccc}} = T^{x_{\sigma \ccc}}$ and $T^{(sx)_{[c,sc]}} = T^{x_{[\sigma c, \sigma sc]}}$.  If $\sigma = \sigma_k$ and not both of $k, k+1 \in \ccc$, then the product $\prod_{c \in \ccc^-}T^{x_{[c,sc]}}$  is exactly the product $\prod_{d \in \sigma(\ccc)^-} T^{x_{[\sigma d, \sigma sd]}}$.  Otherwise, the two products differ by a factor of $T^{s'(x_{[k,k]})} = 1_{A^{s'(x_{[k,k]})}} = 1_I$, and hence can be ignored.  This shows $\dom(sx,\ccc)=\dom(x,\sigma \ccc)$.  We also have $\cod(sx,\ccc) = T^{sx} = T^{x} = \cod(x, \sigma \ccc)$ by assumption.

\item This last item follows from item (iii) of \ref{part1} and simply composing face maps.

\end{enumerate}

\end{proof}

\medskip

We are now ready to state and prove the converse to proposition \ref{forward}.

\medskip

\begin{proposition}\label{backward}

Let $X$ be a simplicial set, and suppose we are given: 

\begin{enumerate}

\item A category $A^x$ for each $1$-simplex $x \in \ds_1$  such that: 

\begin{description}

\item[\textnormal{\textit{(1.a) }}] $A^{sy} = I$ for the degeneracy of a $0$-simplex $y \in \ds_0$.

\end{description}

\item A functor \[T^x: \prod \limits_{i \in \n^-} A^{x_{i,i+1}} \to A^{x_{0,n}}\] 
\noindent for each $n \ge 1$ and $x \in \ds_n$ such that: 

\begin{description}

\item[\textnormal{\textit{(2.a) }}] $T^x = 1_{A^x}$ when $n = 1$.

\item[\textnormal{\textit{(2.b) }}] $T^{sx} = T^{x}$ for any degeneracy map $s$.

\end{description}

\item A natural transformation 

\[\eta^x_{\ccc}: T^{x_\ccc}\circ \prod \limits_{c\in\ccc^-}T^{x_{[c,sc]}} \tto T^x  \] 

\noindent for each $n \ge 1$,   $x \in \ds_n$, and $\{0,n\} \subseteq \ccc \subseteq \n $ such that:

\begin{description}

\item[\textnormal{\textit{(3.a) }}] $\eta^x_{\ccc} = 1_{T^x}$ whenever $\ccc = \{0,n\}$ or $\ccc = \n$. 

\item[\textnormal{\textit{(3.b) }}] $\eta^{sx}_{\ccc} = \eta^{x}_{\sigma \ccc}$ for any degeneracy map $s$.

\end{description}

\end{enumerate}

Then, there exists a unique map $\phi: \ds \to \hcn$ extending the above data with...

\begin{enumerate}
\item $\phi_1(x)(\{0,1\}) = A^x$ for each $x \in \ds_1$.
\item $\phi_n(x)(\{0,n\}\subseteq \n) = T^x$ for each $n \ge 1$ and $x \in \ds_n$.
\item $\phi_n(x)(\{0,n\}\subseteq \ccc \subseteq \n) = \eta^x_{\ccc}$ for each $n \ge 1$, $x \in \ds_n$, and $\lbrace {0,n} \rbrace \subseteq \ccc \subseteq \n$.
\end{enumerate}

\medskip

...\textbf{if and only if} for every $n \ge 4$ and non-degenerate $x\in \ds_n$ and non-degenerate $\left(\lbrace 0,n \rbrace \subset \aaa \subset \bbb \subset \n\right) \in S[n](0,n)_3$, we have:

\begin{equation}\label{key3}\tag*{$(\dagger)$}
\eta^x_{\bbb}\bullet \left(\eta^{x_{\bbb}}_{\delta_{\bbb}^{-1}\aaa} \circ 1 \right)= \eta^x_{\aaa} \bullet \left( 1\circ \prod \limits_{a \in \aaa^-} \eta^{x_{[a,sa]}}_{\delta_{[a,sa]}^{-1}\bbb\cap[a,sa]} \right)     
\end{equation}

\medskip
 Where each $\eta$ fits into one of the four triangular faces of the diagram below \footnote{See in particular item (iii) of lemma \ref{part2}}: 
\medskip

\begin{center}\begin{tikzpicture}[node distance=3.5cm]\label{theorema}
\node (A1) {$\prod \limits_{i \in \n^-}A^{x_{i,i+1}}$};
\node (A2) [right of=A1] {$ $};
\node (A3) [right of=A2] {$ \prod \limits_{a \in \aaa^-} A^{x_{a,sa}}$};
\node (B1) [below of=A1] {$ \prod \limits_{b \in \bbb^-} A^{x_{b,sb}}$};
\node (B2) [right of=B1] {$ $};
\node (B3) [right of=B2] {$A^{x_{0,n}}$};

 \draw[->,  above] (A1) to node {$\prod\limits_{a\in \aaa^-} T^{x_{[a,sa]}}$} (A3);
 \draw[->,  left] (A1) to node {$ \prod\limits_{b \in \bbb^-} T^{x_{[b,sb]}}$} (B1);
 \draw[->,  above, bend left=20] (A1) to node {$ T^{x}$} (B3);
 \draw[->,  below] (B1) to node {$T^{x_\bbb}$} (B3);
 \draw[->,  right] (A3) to node {$ T^{x_\aaa} $} (B3);
 \draw[->,  left] (B1) to node {$\prod\limits_{a \in \aaa^-} T^{x_{\bbb \cap [a,sa]}}$} (A3);
   
\end{tikzpicture}\end{center}

\end{proposition}

\medskip

\begin{proof}

Given an extension $\phi$ (which is necessarily unique by proposition \ref{forward}) the equation \ref{key3} precisely asserts that $\phi_n(x)(\lbrace 0,n \rbrace \subset \aaa \subset \bbb \subset \n)$ is indeed a $3$-simplex of $N_2(\mathsf{Cat})(*,*)$.  That is, it asserts that the diagram of natural transformations associated to the boundary of $\phi_n(x)(\lbrace 0,n \rbrace \subset \aaa \subset \bbb \subset \n)$ is commutative.

Conversely, given the data of items (1), (2), and (3), such that equation \ref{key3} is true, we must assign to each $x \in \ds_n$ an enriched functor $\phi_n(x):S[n]\to N_2(\mathsf{Cat})$.  By lemma \ref{Llemma}, we must therefore give assignments:

\begin{enumerate}

\item A category $\phi_n(x)(\lbrace p, q \rbrace)$ for all $0 \le p < q \le n$
\item A functor 
\[ \phi_n(x)(\lbrace p,q \rbrace \subseteq \ccc):  \prod \limits_{c\in \ccc^-} \phi_n(x)(\lbrace c,sc \rbrace) \to \phi_n(x)(\lbrace p,q \rbrace)\] for each $ \lbrace p,q \rbrace \subseteq \ccc$.
\item A transformation \[\phi_n(x)(\lbrace p,q \rbrace \subseteq \ccc_1 \subseteq \ccc_2): \phi_n(x)(\lbrace p,q \rbrace \subseteq \ccc_1) \circ \prod \limits_{c \in \ccc_1^-} \phi_n(x)(\lbrace c,sc \rbrace \subseteq \ccc_2 \cap [c,sc]) \tto \phi_n(x)(\lbrace p,q \rbrace \subseteq \ccc_2)\]for each  $\lbrace p,q \rbrace \subseteq \ccc_1 \subseteq \ccc_2$.

\end{enumerate}

\noindent Then for an arbitrary $x \in \ds_n$, we define:

\begin{enumerate}

\item $ \phi_n(x)(\{p,q\}):=A^{x_{p,q}}$

\item $\phi_n(x)(\{p,q\}\subseteq \ccc):= T^{x_{\ccc}}$

\item $\phi_n(x)(\{p,q\}\subseteq \ccc_1 \subseteq \ccc_2):= \eta^{x_{\ccc_2}}_{\delta_{\ccc_2}^{-1}\ccc_1}$

\end{enumerate}

Note first that $\phi_n(x)(\{p,q\}\subseteq \ccc)=T^{x_{\ccc}}$ should be a functor $\prod_{c\in \ccc^-} \phi_n(x)(\lbrace c,sc \rbrace) \to \phi_n(x)(\lbrace p,q \rbrace)$, or rather, should be a functor:

\[T^{x_{\ccc}}:\prod_{c\in \ccc^-} A^{x_{c,sc}} \to A^{x_{p,q}}\]  

\noindent This is so by item (iii) of lemma \ref{part1}, as the stipulations $(2.a)$ and $(2.b)$ of the statement of the proposition verify the hypotheses of that lemma.  The natural transformation $\eta^{x_{\ccc_2}}_{\delta_{\ccc_2}^{-1}\ccc_1}$ should be a natural transformation with domain and codomain:

\[\eta^{x_{\ccc_2}}_{\delta_{\ccc_2}^{-1}\ccc_1}: T^{x_{\ccc_1}} \circ \prod\limits_{c'\in \ccc_1}T^{x_{\ccc_2 \cap [c',sc']}} \tto T^{x_{\ccc_2}} \]

\noindent This follows from item (iii) of lemma \ref{part2}, as again the additional stipulations $(3.a)$ and $(3.b)$ verify the hypotheses of that lemma.

Note also that to use lemma \ref{Llemma}, we must have that $\phi_n(x)(\lbrace p,q \rbrace = \lbrace p,q \rbrace) = 1_{\phi_n(x)(\lbrace p,q \rbrace)}$.  This follows from the stipulation $(2.a)$.  Similarly, lemma \ref{Llemma} also requires that $\phi_n(x)(\lbrace p,q \rbrace \subseteq \ccc_1 \subseteq \ccc_2)$ is $1_{\phi_n(x)(\lbrace p,q \rbrace \subseteq \ccc_2)}$  if either `$\subseteq$' is an `$=$'.  This follows from the stipulation $(3.a)$.

By lemma \ref{Llemma}, these assignments extend to an enriched functor $\phi_n(x)$ if for every $(\{p,q\}\subset \ccc_1 \subset \ccc_2 \subset \ccc_3) \in S[n](0,n)_3$, we have:

\begin{align*}\label{key4} & \phi_n(x)(\{p,q\} \subset \ccc_2 \subset \ccc_3) \bullet (\phi_n(x)( \{p,q\} \subset \ccc_1 \subset \ccc_2) \circ 1) = \\& \phi_n(x)(\{p,q\} \subset \ccc_1 \subset \ccc_3) \bullet \left( 1 \circ \prod \limits_{c \in \ccc_1^-} \phi_n(x)(\{c,sc\} \subseteq \ccc_2\cap [c,sc] \subseteq \ccc_3 \cap [c,sc])\right)
\end{align*}

\noindent Tracing through the definitions and the assignments in lemma \ref{Llemma}, this is precisely the equation \ref{key3} with respect to $\delta_{\ccc_3}^{-1}(\{p,q\}\subset \ccc_1 \subset \ccc_2 \subset \ccc_3)$, and is hence true by hypothesis.

We have thus far shown that to each $x \in \ds_n$, we can assign an enriched functor $\phi_n(x) \in \hcn_n$. We must now check that $\phi$ commutes with face and degeneracy maps and so is a simplicial map.  For each $x\in \ds_{n+1}$ and $\delta:\n\to [n+1]$ with corresponding face map $d$, we must show that $\phi_{n}(dx) = \phi_{n+1}(x) \circ \delta:S[n]\to N_2(\mathsf{Cat})$.  As such, it suffices to check that their action is the same on generating simplices.  This follows more or less by definition, and we show only the $2$-dimensional generating simplex case:

\begin{align*}
\phi_n(dx)(\{p,q\} \subseteq \ccc_1 \subseteq \ccc_2) =& \eta^{dx_{\ccc_2}}_{\delta_{\ccc_2}^{-1}\ccc_1}  \\
 =& \eta^{x_{\delta \ccc_2}}_{\delta_{\delta \ccc_2}^{-1} \delta \ccc_1}\\
 =& \phi_{n+1}(x)(\{\delta p, \delta q\} \subseteq \delta \ccc_1 \subseteq \delta \ccc_2 )
\end{align*}

\medskip

For degeneracy maps we must show that for arbitrary $x \in \ds_{n-1}$ and $\sigma:\n\to [n-1]$, we have $\phi_{n}(sx) = \phi_{n-1}(x) \circ \sigma:S[n]\to N_2(\mathsf{Cat})$.  Again it suffices to check equality on generating simplices. Recalling lemma \ref{degens}, we have:

\begin{enumerate}[(i)]
\item $\begin{aligned}[t]
    \phi_{n}(sx)(\{p,q\}) &= A^{(sx)_{p,q}}=A^{s'(x_{\sigma p, \sigma q}) } = A^{x_{\sigma p, \sigma q} }= \left(\phi_{n-1}(x) \circ \sigma \right) (\{p,q\}) \\
\end{aligned}$\\
\item $\begin{aligned}[t]
    \phi_{n}(sx)(\{p,q\} \subseteq \ccc) &= T^{(sx)_{\ccc}}=T^{s'(x_{\sigma\ccc})} = T^{x_{\sigma \ccc} }= \left( \phi_{n-1}(x) \circ \sigma \right)(\{p,q\}\subseteq \ccc) \\
\end{aligned}$\\
\item $\begin{aligned}[t]
    \phi_{n}(sx)(\{p,q\} \subseteq \ccc_1 \subseteq \ccc_2) &= \eta^{(sx)_{\ccc_2}}_{\delta_{\ccc_2}^{-1}\ccc_1} =\eta^{s'(x_{\sigma \ccc_2})}_{\delta_{\ccc_2}^{-1}\ccc_1}\\ 
    &= \eta^{x_{\sigma \ccc_2}}_{ \sigma' \delta_{\ccc_2}^{-1}\ccc_1} = \eta^{x_{\sigma \ccc_2}}_{\delta_{\sigma\ccc_2}^{-1}\sigma\ccc_1} = \left(\phi_{n-1}(x) \circ \sigma \right)(\{p,q\}\subset \ccc_1 \subset \ccc_2)\\
\end{aligned}$
\end{enumerate}

The fourth equation of item (iii) follows because we have $\delta_{\sigma \ccc_2}  \sigma' = \sigma \delta_{\ccc_2}$ as in the proof of lemma \ref{degens}.

\end{proof}

\medskip

We will make use of the general version of proposition $\ref{backward}$, but care most about the specific case when $X = \ctn$.  In this case the proposition simplifies nicely because there are only two $1$-simplices, $0,1 \in \ctn_1$, and only $1$ is non-degenerate. Furthermore, in light of the results to follow, we can think an arbitrary map $\phi:\ctn \to \hcn$ as a truly general kind of monoidal-type category.  When $X=\ctn$, proposition \ref{backward} can be taken to be a presentation of a monoidal-type category in the usual dichotomy of data  (1), (2), and (3),  subject to coherence  \ref{key3}.  Writing $\sum \spi{x}$ for the total number of `$1$'s in $\spi{x}$ for a simplex $x \in \ctn$, we have that an arbitrary map $\phi: \ctn \to \hcn$ can be defined by:

\begin{enumerate}
\item A category $A = A^1$, and write $A^0=I$.

\item A functor  

 \[T^x:\prod \limits_{i \in \n^-} A^{x_{i,i+1}} = A^{\sum \spi{x}} \to A\]

for each $n \ge 1$ and $x \in \ctn_n$ subject to $(2.a)$ and $(2.b)$.

\item A natural transformation

\[ \eta^x_{\ccc} : (T^{x_\ccc} \circ \prod\limits_{c \in \ccc^-}T^{x_{[c,sc]}}) \tto T^x \] 

for each face $n \ge 1$,  $x \in \ctn_n$, and $\lbrace 0,n \rbrace \subseteq \ccc \subseteq \n$ subject to $(3.a)$ and $(3.b)$.

\end{enumerate}

\noindent The natural transformations of item $(3)$ must be \textit{coherent} insofar as the equation \ref{key3} is satisfied for every $\lbrace 0,n \rbrace \subset \aaa \subset \bbb \subset \n$:

\[\eta^x_{\bbb}\bullet \left(\eta^{x_{\bbb}}_{\delta_{\bbb}^{-1}\aaa} \circ 1 \right)= \eta^x_{\aaa} \bullet \left( 1\circ \prod \limits_{a \in \aaa^-} \eta^{x_{[a,sa]}}_{\delta_{[a,sa]}^{-1}\bbb\cap[a,sa]} \right)     \]

 \medskip

\section{Classifying lax monoidal categories}

In this section, we will prove proposition \ref{laxclass}.  Given a lax monoidal category $(A,\otimes^n, \iota, \gamma)$, we first define the classifying map $\alpha$, then verify its well definedness, and lastly show that given $\alpha$ we can reconstruct $(A,\otimes^n, \iota, \gamma)$.

\medskip

\begin{definition}\textbf{The Lax Monoidal Category Classifying Map}

\medskip

Let $(A,\otimes^n, \iota, \gamma)$ be a lax monoidal category.  We assign to $A$ the map $\alpha: \ctn \to \hcn $ using proposition \ref{backward}.  We define the map as follows:

\begin{enumerate}
\item  Let $A^x := A$ if $x=1$ and $I$ otherwise.  

\medskip

\item Let $T^x := \otimes^{\sum \spi{x}}: A^{\sum \spi{x}} \to A$ if $x$ is non-degenerate with dimension $\ge 2$.  Otherwise $T^x$ is defined in accord with stipulations $(2.a)$ and $(2.b)$ of proposition $\ref{backward}$.  In particular, $T^x := 1_{A^x}$ if $x$ has dimension $1$, or is the degeneracy of a $1$-simplex.

\medskip

\item  For each non-degenerate $x \in \ctn_n$ with $n \ge 3$ and $\lbrace 0,n \rbrace \subset \ccc \subset \n$, let \[\eta^x_{\ccc}:= \gamma \bullet (\overline{\iota} \circ \prod \overline{\iota}): (T^{x_\ccc} \circ \prod\limits_{c \in \ccc^-}T^{x_{[c,sc]}}) \tto T^x\]

\noindent Where $\gamma$ is short-hand for the appropriately sized associativity transformation $\gamma_{n,k_1,...,k_n}$, and $\overline{\iota}:T^x \to \overline{T^x}$ is defined by:

\begin{displaymath}
    \overline{T^{x}} := \left\{
     \begin{array}{lr}
       \otimes^1   &\text{if } T^{x} = 1_A
       \vspace{.2cm}
       \\
	  T^{x}  &\text{else} 
     \end{array}
    \right.  \hspace{2cm}\overline{\iota}:T^{x} \tto \overline{T^{x}} := \left\{
     \begin{array}{lr}
       \iota   &\text{if } T^{x} = 1_A
       \vspace{.2cm}
       \\
	  1_{T^{x}}  &\text{else} 
     \end{array}
   \right.
\end{displaymath} 

\noindent Otherwise $\eta^x_{\ccc}$ is defined in accord with stipulations $(3.a)$ and $(3.b)$ of proposition \ref{backward}.  In particular, we may have $\eta^x_{\ccc} = 1_{T^x}$ if, for example, $x \in \ctn_n$ and $\ccc= \n$ or $\ccc = \{0,n\}$. 

\end{enumerate}

This transformation $\gamma \bullet (\overline{\iota} \circ \prod \overline{\iota})$ encodes a simple idea: $\iiota \circ \prod \iiota$ converts all $1_A$'s appearing in any factor of either composite of the domain into $\otimes^1$'s, so that by the time $\gamma$ is applied, only $\otimes^n$ functor factors remain. 
\end{definition}

\medskip

This data extends to a map $\alpha: \ctn \to \hcn$ via proposition \ref{backward} because we have the following:

\begin{proposition}\label{lemma1}
Let $(A,\otimes^n, \iota, \gamma)$ be a lax monoidal category, and let $A^x, T^x, \eta^x_{\ccc}$ be defined as above.   Then for every $n \ge 4$,  $\{0,n\} \subset \aaa \subset \bbb \subset \n$ and non-degenerate $x \in \ctn_n$, equation \ref{key3} is satisfied:  

\[ \eta^x_{\bbb}\bullet \left(\eta^{x_{\bbb}}_{\delta_{\bbb}^{-1}\aaa} \circ 1 \right)= \eta^x_{\aaa} \bullet \left( 1\circ \prod \limits_{a \in \aaa^-} \eta^{x_{[a,sa]}}_{\delta_{[a,sa]}^{-1}\bbb\cap[a,sa]} \right)      \]

\end{proposition}

\medskip

As $x\in \ctn_n$ is a fixed non-degenerate simplex throughout the proof, we will simply omit `$x$' from our notation.  We will write $A_{i,i+1}:= A^{x_{i,i+1}}$,
$T:= T^x$, $T_{\ccc}:=T^{x_{\ccc}}$ and $\eta_{\ccc}:= \eta^x_{\ccc}$  for any $\{0,n\} \subseteq \ccc \subseteq \n $.  We write $\eta^{\ccc}_{\ddd}:= \eta^{x_{\ccc}}_{\delta_{\ccc}^{-1}\ddd}$ for any $\{0,n\} \subseteq \ddd \subseteq \ccc \subseteq [n]$.  In this notation, we will prove:

\[\eta_{\bbb} \bullet (\eta^{\bbb}_{\aaa} \circ 1) = \eta_{\aaa} \bullet \left( 1 \circ \prod \limits_{a\in \aaa^-} \eta^{[a,sa]}_{\bbb \cap [a,sa]} \right) \]

 \noindent It will be convenient to recall that these natural transformations  each fit into one of the four triangular faces of:

\begin{center}\begin{tikzpicture}[node distance=3.5cm]\label{theorem}
\node (A1) {$\prod \limits_{i \in \n^-}A_{i,i+1}$};
\node (A2) [right of=A1] {$ $};
\node (A3) [right of=A2] {$ \prod \limits_{a \in \aaa^-} A_{a,sa}$};
\node (B1) [below of=A1] {$\prod \limits_{b \in \bbb^-} A_{b,sb}=   \hspace{-.35cm} \prod \limits_{\substack{ a\in \aaa^- \\ b \in \bbb \cap [a,sa]^-}} \hspace{-.35cm} A_{b,sb}$};
\node (B2) [right of=B1] {$ $};
\node (B3) [right of=B2] {$A_{0,n}$};

 \draw[->,  above] (A1) to node {$\prod\limits_{a\in \aaa^-} T_{[a,sa]}$} (A3);
 \draw[->,  left] (A1) to node {$ \prod\limits_{b \in \bbb^-} T_{[b,sb]}$} (B1);
 \draw[->,  above, bend left=20] (A1) to node {$ T$} (B3);
 \draw[->,  below] (B1) to node {$T_{\bbb}$} (B3);
 \draw[->,  right] (A3) to node {$ T_{\aaa} $} (B3);
 \draw[->,  left] (B1) to node {$\prod\limits_{a \in \aaa^-} T_{\bbb \cap [a,sa]}$} (A3);
   
\end{tikzpicture}\end{center}

\begin{center} \hspace{.4cm} Figure $(\ddagger)$ \end{center}

\medskip

\begin{proof} (of proposition \ref{lemma1})

Because $x$ is non-degenerate and there are no `$=$'s in $\{0,n\} \subset \aaa \subset \bbb \subset \n$ , we know that both $\eta_{\aaa}$ and $\eta_{\bbb}$ are of the form $\gamma \bullet \iiota$ by definition.  However, $x_{\bbb}$ or $x_{[a,sa]}$ may be degenerate, and/or we may have `$=$' in $\{a,sa\} \subseteq \bbb \cap [a,sa] \subseteq [a,sa]$ for any $a \in \aaa^-$.  As a result, $\eta^{\bbb}_{\aaa}$ or $\eta^{[a,sa]}_{[a,sa]\cap \bbb}$ may be given via stipulations $(3.a)$ and/or $(3.b)$ as identity natural transformations.  For example, if $x_{\bbb}$ was the degeneracy of a $2$-simplex $y$, then  $\eta^{\bbb}_{\aaa} = 1_{T^y}$ as can easily be checked.  This shows we will need to think about how all three of $\gamma$, $\iota$, and identity natural transformations interact in the proof.  Luckily, we need only consider two cases: when $\eta^{\bbb}_{\aaa} = \gamma \bullet \iiota$, and when it is an identity.  

Let us first assume that $\eta^{\bbb}_{\aaa} = \gamma \bullet \iiota$.  Define the set:

\[ \Gamma = \lbrace a \in \aaa^- ~|~ \eta^{[a,sa]}_{\bbb \cap [a,sa]} \text{ is defined by } \gamma \bullet \iiota \rbrace \] 

\noindent So that $a \in (\aaa^- )$\textbackslash$\Gamma$ only if $ \eta^{[a,sa]}_{\bbb \cap [a,sa]} = 1_{T_{[a,sa]}}$.  Consider the following diagram:
\\\\\\\\\\\\\\\\\\\\\\\\\\\\\\\\\\\\\\\\\\\\\\\\\\\\\\\\\\\\\\\\\\\\\\\\\\\\\\\\\\\\\\\\\\\

 \pagebreak

\vspace{-2cm}

\begin{center}\begin{sideways}\vspace{-7cm}\begin{tikzpicture}[node distance=3cm]\vspace{-8cm}\small

  \node(A0) {$ \substack{T_{\aaa} \circ \prod \limits_{a\in \aaa^-} T_{\bbb \cap [a,sa]} \circ \prod \limits_{b\in \bbb^-} T_{[b,sb]} =\\  T_{\aaa} \circ \prod \limits_{a\in \aaa^-} \left( T_{\bbb \cap [a,sa]} \circ \prod \limits_{b \in (\bbb \cap [a,sa])^-} T_{[b,sb]} \right) } $} ;  
  \node (A1) [right of=A0] {$ $} ;
  \node (A2) [right of=A1] {$ $} ;
  \node (A3) [right of=A2] {$ $} ;
  \node (A4) [right of=A3] {$ T_{\aaa} \circ \left( \prod \limits_{a \in \Gamma} \left( T_{\bbb \cap [a,sa]} \circ \prod \limits_{b \in (\bbb \cap [a,sa])^-} T_{[b,sb]} \right) \times  \prod \limits_{a \nin \Gamma} \left( T_{\bbb \cap [a,sa]} \circ \prod \limits_{b \in (\bbb \cap [a,sa])^-} T_{[b,sb]} \right) \right)$} ;
  
  \node (B0) [below of=A0]{$T_{\aaa} \circ \left( \prod \limits_{a\in \Gamma} T_{\bbb \cap [a,sa]} \times \prod \limits_{a \nin \Gamma}  T_{\bbb \cap [a,sa]} \right) \circ \left( \prod \limits_{\substack{b \in (\bbb \cap [a,sa])^-\\ a \in \Gamma}} \hspace{-.2cm} T_{[b,sb]} \times \hspace{-.1cm} \prod \limits_{\substack{b \in (\bbb \cap [a,sa])^- \\ a \nin \Gamma}} \hspace{-.2cm} T_{[b,sb]}\right) $} ;  
  \node (B1) [right of=B0] {$   $} ;
  \node (B2) [right of=B1] {$ $} ;
  \node (B3) [right of=B2] {$ $} ;
  \node (B4) [right of=B3] {$ \overline{T_{\aaa}} \circ \left( \prod \limits_{a \in \Gamma} \left( \overline{T_{\bbb \cap [a,sa]}} \circ \prod \limits_{b \in (\bbb \cap [a,sa])^-} \overline{T_{[b,sb]}} \right) \times  \prod \limits_{a \nin \Gamma} \left( T_{\bbb \cap [a,sa]} \circ \prod \limits_{b \in (\bbb \cap [a,sa])^-} T_{[b,sb]} \right) \right) $} ;
 
 \node (E0) [below of=B0]{$ $} ;  
  \node (E1) [right of=E0] {$ $} ;
  \node (E2) [right of=E1] {$\overline{T_{\aaa}} \circ \left( \prod \limits_{a \in \Gamma} \left( \overline{T_{\bbb \cap [a,sa]}} \circ \prod \limits_{b \in (\bbb \cap [a,sa])^-} \overline{T_{[b,sb]}} \right) \times  \prod \limits_{a \nin \Gamma} \left( \overline{T_{\bbb \cap [a,sa]}} \circ \prod \limits_{b \in (\bbb \cap [a,sa])^-} \overline{T_{[b,sb]}} \right) \right) $} ;
  \node (E3) [right of=E2] {$ $} ;
  \node (E4) [right of=E3] {$ $} ;

  \node (C0) [below of=E0]{$ T_{\bbb}\circ \left(\prod \limits_{\substack{b \in (\bbb \cap [a,sa])^-\\ a \in \Gamma}} T_{[b,sb]}\times \prod \limits_{\substack{b \in (\bbb \cap [a,sa])^-\\ a \nin \Gamma}} T_{[b,sb]}\right)$} ;  
  \node (C5) [below of=E0, below]{$ $};
  \node (C1) [right of=C0] {$ $} ;
  \node (C2) [right of=C1] {$ $} ;
  \node (C3) [right of=C2] {$ $} ;
  \node (C4) [right of=C3] {$  \overline{T_{\aaa}} \circ \left( \prod \limits_{a \in \Gamma} \overline{T_{[a,sa]}}\times \prod \limits_{a \nin \Gamma} T_{[a,sa]}  \right)$} ;

  \node (D0) [below of=C0]{$ T_{\bbb}\circ \left(\prod \limits_{\substack{b \in (\bbb \cap [a,sa])^-\\ a \in \Gamma}} \overline{T_{[b,sb]}}\times \prod \limits_{\substack{b \in (\bbb \cap [a,sa])^-\\ a \nin \Gamma}} \overline{T_{[b,sb]}}\right)$} ;  
  \node (D1) [right of=D0] {$ $} ;
  \node (D2) [right of=D1] {$T$} ;
  \node (D3) [right of=D2] {$ $} ;
  \node (D4) [right of=D3] {$ \overline{T_{\aaa}} \circ \left( \prod \limits_{a \in \Gamma} \overline{T_{[a,sa]}}\times \prod \limits_{a \nin \Gamma} \overline{T_{[a,sa]}}    \right)$} ;

  \node (F0) [below of=D0]{$ $} ;  
  \node (F1) [right of=F0] {$ $} ;
  \node (F2) [right of=F1] {$ $} ;
  \node (F3) [right of=F2] {$ $} ;
  \node (F4) [right of=F3] {$ $} ;

  \draw[=, white](C2) to node [black, above] {$ \textbf{S} $} (D2);
  \draw[=, white](B2) to node [black] {$ \textbf{N} $} (A2);
  \draw[=, white](C2) to node [black] {$ \textbf{W} $} (E0);
  \draw[=, white](C2) to node [black] {$ \textbf{E} $} (E4);
 \draw[=,line width=.06cm](A0) to node [above] {$1$} (A4);
 \draw[=, white, line width=.03cm](A0) to node [below] {$ $} (A4);
 \draw[=, line width=.06cm](A0) to node [left] {$1$} (B0);
 \draw[=, white, line width=.03cm](A0) to node [above] {$  $} (B0);
 \draw[->, bend right=15](B0) to node [left] {$ \eta^{\bbb}_{\aaa} \circ 1$} (C0);
 \draw[->](B0) to node [above right] {$ \overline{\iota} \circ \left(\overline{\iota} \times \overline{\iota} \right) \circ \left(\overline{\iota} \times \overline{\iota} \right) $} (E2);
 \draw[->](A4) to node [right] {$ \overline{\iota} \circ \left(\left(\overline{\iota} \circ \overline{\iota} \right) \times \left(1 \circ 1 \right)\right) $} (B4);
 \draw[->](B4) to node [above left] {$ 1\circ \left(\left(1 \circ 1 \right) \times \left(\overline{\iota} \circ \overline{\iota} \right)\right) $} (E2);
 \draw[->, bend left=15](B4) to node [right] {$1 \circ \left( \prod \gamma \times 1\right) $} (C4); 
 \draw[->](C4) to node [right] {$ 1 \circ \left( 1\times \overline{\iota} \right)$} (D4);  
 \draw[->](D4) to node [below] {$\gamma$} (D2); 
 \draw[->](C0) to node [left] {$ 1 \circ \left( \overline{\iota} \times \overline{\iota}\right) $} (D0); 
 \draw[->](D0) to node [below] {$\gamma$} (D2);  
 \draw[->, bend right=30](E2) to node [right] {$ 1 \circ \left( \prod \gamma \times \prod \gamma \right)$} (D4); 
 \draw[->, bend left=35](E2) to node [left] {$ \gamma \circ \left(1 \times 1 \right) $} (D0);

\end{tikzpicture}\end{sideways}\end{center}

\pagebreak

Due to interchange of horizontal composition with products, the middle node in this diagram can be rewritten:

\[\overline{T_{\aaa}} \circ \left( \prod \limits_{a \in \Gamma} \left( \overline{T_{\bbb \cap [a,sa]}} \circ \prod \limits_{b \in (\bbb \cap [a,sa])^-} \overline{T_{[b,sb]}} \right) \times  \prod \limits_{a \nin \Gamma} \left( \overline{T_{\bbb \cap [a,sa]}} \circ \prod \limits_{b \in (\bbb \cap [a,sa])^-} \overline{T_{[b,sb]}} \right) \right) =\] 

\[  \overline{T_{\aaa}} \circ \left( \prod \limits_{a\in \Gamma} \overline{T_{\bbb \cap [a,sa]}} \times \prod \limits_{a \nin \Gamma}  \overline{T_{\bbb \cap [a,sa]}} \right) \circ \left( \prod \limits_{\substack{b \in (\bbb \cap [a,sa])^- \\ a \in \Gamma}}  \overline{T_{[b,sb]}} \times \hspace{-.1cm} \prod \limits_{\substack{b \in (\bbb \cap [a,sa])^- \\ a \nin \Gamma}}  \overline{T_{[b,sb]}}\right)\]
\\
\noindent The map out of this central node $\gamma \circ (1 \times 1)$ is written in regards to this second expression.  The notation for the other map out of this central node, $1 \circ \prod \gamma \times \prod \gamma$ is written with regards to the first.
\\
\\
First we claim that commutativity of the above diagram implies proposition \ref{theorem}.  As we have assumed that $\eta^{\bbb}_{\aaa}$ is of the form $\gamma \bullet \iiota$, we have that $x_{\bbb}$ is not the degeneracy of a $2$-simplex, and so not the degeneracy of a $1$-simplex, and hence that $T_{\bbb}$ is not an identity map.  Thus $T_{\bbb}=\overline{T_{\bbb}}$.  The composition of the left most edges is then:

\[ \gamma \bullet \left(1 \circ \left( \overline{\iota} \times \overline{\iota}\right) \right) \bullet  \left(\eta^{\bbb}_{\aaa} \circ 1\right) = \eta_{\bbb} \bullet  \left( \eta^{\bbb}_{\aaa} \circ 1\right)\]
\\

The composition of the right most edges is:

\begin{align}
&\gamma \bullet \left( 1 \circ \left( 1\times \overline{\iota} \right)\right) \bullet \left( 1 \circ \left( \prod \gamma \times 1\right) \right)  \bullet \left(\overline{\iota} \circ \left(\left(\overline{\iota} \circ \overline{\iota} \right) \times \left(1 \circ 1 \right)\right) \right)\\    \label{1}
=&~\gamma \bullet \left( \overline{\iota} \circ \left( 1\times \overline{\iota} \right) \right)\bullet \left( 1 \circ \left( \prod \gamma \times 1\right)  \right) \bullet \left( 1 \circ \left(\left(\overline{\iota} \circ \overline{\iota} \right) \times \left(1 \circ 1 \right)\right)\right) \\ 
\label{2}
=&~\gamma \bullet \left( \overline{\iota} \circ \left( \overline{\iota} \times \overline{\iota} \right) \right) \bullet \left( 1 \circ \left( \prod \gamma \times 1\right) \right)  \bullet \left( 1 \circ \left(\left(\overline{\iota} \circ \overline{\iota} \right) \times \left(1 \circ 1 \right)\right)\right) \\
\label{3}
=&~\eta_{\aaa} \bullet \left( 1 \circ \left( \prod \gamma \times 1\right)\right)  \bullet \left( 1 \circ \left(\left(\overline{\iota} \circ \overline{\iota} \right) \times \left(1 \circ 1 \right)\right) \right)  \\
\label{4}
=&~\eta_{\aaa} \bullet \left( 1 \circ \left(\left( \prod \gamma \bullet \left( \overline{\iota} \circ \overline{\iota} \right) \right) \times \left( 1 \bullet \left( 1 \circ 1 \right)\right) \right)\right) \\
\label{5}
=&~\eta_{\aaa} \bullet \left( 1 \circ \left(\left( \prod_{a \in \Gamma} \eta \right) \times \left( \prod_{a \nin \Gamma} \eta \right) \right) \right) \\
\label{6}
=&~\eta_{\aaa} \bullet \left( 1 \circ \prod \limits_{a\in \aaa^-} \eta^{[a,sa]}_{\bbb \cap [a,sa]} \right) 
		\end{align}

\medskip

\noindent Equation \ref{1} follows from commuting $\overline{\iota}$ with the identity functor $1$.  Equation \ref{2} follows because for each $a\in \Gamma$, the map...

\[\gamma: \left( \overline{T_{\bbb \cap [a,sa]}} \circ \prod \limits_{b \in (\bbb \cap [a,sa])^-} \overline{T_{[b,sb]}} \right)  \tto T_{[a,sa]} \]
\\
... cannot possibly have as codomain $T_{[a,sa]} = 1_A $, as the output of $\gamma$ is always a tensor $\otimes$.  Hence we have $T_{[a,sa]} = \overline{T_{[a,sa]}}$ for every $a \in \Gamma$.  This means:

\[1 = \overline{\iota} : T_{[a,sa]} \tto \overline{T_{[a,sa]}}\]
\\
Equation \ref{3} is then the assumption that $\eta_{\aaa}$ is of the form $\gamma \bullet \iiota \circ \prod \iiota$.  Equation \ref{4} is repeated applications of interchange of vertical composition.  Equation \ref{5} is just the definition of $\eta$ and $\Gamma$:  The first grouping corresponds to $\eta$ for $a \in \Gamma$ which in turn corresponds to all those $\eta$ of the form $\gamma \bullet \overline{\iota}$;  The second grouping corresponds to all of the other $\eta$'s, which are thus simply $1$.
\\
\\
Thus the commutativity of the above diagram implies proposition $\ref{theorem}$.  We verify the commutativity by verifying the commutativity of each sub diagram: \textbf{N, W, S, E}. 
\\
\begin{enumerate}[ ] 

\item \textbf{N}:  This commutes trivially.
\\
\\
\item \textbf{W}: This commutes as a result of commuting $\overline{\iota}$ with $1$.

\begin{align} &~\left(\gamma \circ \left(1 \times 1 \right)\right) \bullet \left(\overline{\iota} \circ \left(\overline{\iota} \times \overline{\iota} \right) \circ \left(\overline{\iota} \times \overline{\iota} \right)\right) \\
=&~ \left(\gamma \circ \left(\overline{\iota} \times \overline{\iota} \right)\right) \bullet \left(\overline{\iota} \circ \left(\overline{\iota} \times \overline{\iota} \right) \circ \left(1 \times 1 \right)\right) \\
=&~ \left( 1 \circ \left(\overline{\iota} \times \overline{\iota} \right) \right) \bullet \left( \gamma \circ \left( 1 \times 1 \right)\right)\bullet \left( \overline{\iota} \circ \left(\overline{\iota} \times \overline{\iota} \right) \circ \left( 1 \times 1 \right)\right) \\
=&~ \left( 1 \circ \left(\overline{\iota} \times \overline{\iota} \right)\right) \bullet \left(\gamma \bullet \left( \overline{\iota} \circ \left(\overline{\iota} \times \overline{\iota} \right)\right)\right) \circ \left((1 \times 1) \bullet (1\times 1)\right)\\
=&~ \left(1 \circ \left(\overline{\iota} \times \overline{\iota} \right)\right) \bullet \left(\eta^{\bbb}_{\aaa} \circ \left( 1 \times 1 \right)\right)
\end{align}

\item \textbf{S}: Note that every vertex in this square consists of only products of $\otimes$'s (and $1_I$'s, which again, are ignored).  The commutativity of this square is simply the associativity axiom for $\gamma$ coming from the definition of a lax monoidal category.
\\
\\
\item \textbf{E}:  The eastern square $E$ commutes if and only if the following square commutes:

\hfill\\

\hspace{-1cm}\begin{tikzpicture}[node distance=3.2cm]\label{reducta}

   \node(A0) {$ \prod \limits_{a \in \Gamma} \left( \overline{T_{\bbb \cap [a,sa]}} \circ \prod \limits_{b \in (\bbb \cap [a,sa])^-} \overline{T_{[b,sb]}}\right) \times \prod \limits_{a \nin \Gamma} \left( T_{\bbb \cap [a,sa]} \circ \prod \limits_{b \in (\bbb \cap [a,sa])^-} T_{[b,sb]} \right)   $} ;
  \node (A1) [right of=A0] {$   $} ;
  \node (B1) [below of=A1] {$ $};
  \node (B0) [left of=B1] {$\prod \limits_{a \in \Gamma} \left( \overline{T_{\bbb \cap [a,sa]}} \circ \prod \limits_{b \in (\bbb \cap [a,sa])^-} \overline{T_{[b,sb]}}\right) \times \prod \limits_{a \nin \Gamma} \left(\overline{T_{\bbb \cap [a,sa]}} \circ \prod \limits_{b \in (\bbb \cap [a,sa])^-} \overline{T_{[b,sb]}} \right)  $};
  \node (A2) [right of=A1] {$ $};
  \node (A3) [right of=A2] {$\prod \limits_{a \in \Gamma} \overline{T_{[a,sa]}} \times  \prod \limits_{a \in \Gamma} T_{[a,sa]} $};
  \node (B3) [below of=A3] {$\prod \limits_{a \in \Gamma} \overline{T_{[a,sa]}} \times  \prod \limits_{a \in \Gamma} \overline{T_{[a,sa]}}  $};
  
  \draw[->,  above] (A0) to node {$ \prod \gamma \times \prod 1$} (A3);  
  \draw[->,  left] (A0) to node {$ \prod 1 \times \prod \overline{\iota}\circ \overline{\iota}$} (B0);
  \draw[->,  below] (B0) to node {$ \prod \gamma \times \prod \gamma$} (B3);
  \draw[->,  left] (A3) to node {$ \prod 1 \times \prod \overline{\iota}$} (B3);
\end{tikzpicture}

\hfil\\The square commutes in the first component of the binary product trivially:

\[  \prod \gamma \bullet \prod 1 = \prod 1 \bullet \prod \gamma \]
\\
It suffices then to show that the square commutes in the second component, and decomposing the product $ \prod \limits_{a \nin \Gamma}$, we see it suffices to show that for each $a \nin \Gamma$ we have:

\begin{center}\begin{tikzpicture}[node distance=3.2cm]\label{reductb}

   \node(A0) {$  T_{\bbb \cap [a,sa]} \circ \prod \limits_{b \in (\bbb \cap [a,sa])^-} T_{[b,sb]} $} ;
  \node (A1) [right of=A0] {$  T_{[a,sa]} $} ;
  \node (B1) [below of=A1] {$ \overline{T_{[a,sa]}} $};
  \node (B0) [left of=B1] {$\overline{ T_{\bbb \cap [a,sa]}} \circ \prod \limits_{b \in (\bbb \cap [a,sa])^-} \overline{ T_{[b,sb]}} $};
  \draw[->,  above] (A0) to node {$ 1 $} (A1);  
  \draw[->,  left] (A0) to node {$\overline{\iota} \circ \overline{\iota}$} (B0);
  \draw[->,  below] (B0) to node {$ \gamma $} (B1);
  \draw[->,  right] (A1) to node {$ \overline{\iota} $} (B1);
\end{tikzpicture}\end{center}

\medskip

\noindent  Now, $T_{[a,sa]}$ may be a tensor $\otimes^m$, $1_A$, or $1_I$.  If $T_{[a,sa]}=\otimes^m$, the map $ T_{\bbb \cap [a,sa]} \circ \prod  T_{[b,sb]}$ may be $1_A \circ \otimes^m$ or $\otimes^m \circ \prod 1_A$. The square in question therefore becomes one of:

\medskip

\begin{center}\begin{tikzpicture}[node distance=1.8cm]\label{reductd2}

   \node(A0) {$ \otimes^m = \otimes^m \circ  \prod 1_A $} ;
  \node (A1) [right of=A0] {$  $} ;
  \node (B1) [below of=A1] {$ \otimes^m $};
  \node (B0) [left of=B1] {$\otimes^m \circ \prod \otimes^1$};
  \draw[->,  above right] (A0) to node {$ 1_{\otimes^m} $} (B1);  
  \draw[->,  left] (A0) to node {$\iota$} (B0);
  \draw[->,  below] (B0) to node {$ \gamma $} (B1);
 
\end{tikzpicture}\begin{tikzpicture}[node distance=1.8cm]\label{reductd3}

   \node(A0) {$ \otimes^m = 1_A \circ \otimes^m $} ;
  \node (A1) [right of=A0] {$  $} ;
  \node (B1) [below of=A1] {$ \otimes^m $};
  \node (B0) [left of=B1] {$\otimes^1 \circ \otimes^m$};
  \draw[->,  above right] (A0) to node {$ 1_{\otimes^m} $} (B1);  
  \draw[->,  left] (A0) to node {$\iota$} (B0);
  \draw[->,  below] (B0) to node {$ \gamma $} (B1);
 
\end{tikzpicture}\end{center}

\noindent These triangles are precisely those appearing in the unitality axiom of the lax monoidal category, hence commute by definition. 

In the case when $T_{[a,sa]}=1_A$, then the square reduces to the following:

\begin{center}\begin{tikzpicture}[node distance=1.8cm]\label{reductc}

   \node(A0) {$ 1_A \circ 1_A $} ;
  \node (A1) [right of=A0] {$ 1_A $} ;
  \node (B1) [below of=A1] {$ \otimes^1 $};
  \node (B0) [left of=B1] {$\otimes^1 \circ \otimes^1$};
  \draw[->,  above] (A0) to node {$ 1 $} (A1);  
  \draw[->,  left] (A0) to node {$\iota \circ \iota$} (B0);
  \draw[->,  below] (B0) to node {$ \gamma $} (B1);
  \draw[->,  right] (A1) to node {$ \iota $} (B1);
\end{tikzpicture}\end{center}

\noindent The unitality axiom in the definition of the lax monoidal category gives in particular that the following triangle commutes:

\begin{center}\begin{tikzpicture}[node distance=1.8cm]\label{reductd}

   \node(A0) {$ \otimes^1 = \otimes^1 \circ 1_A $} ;
  \node (A1) [right of=A0] {$  $} ;
  \node (B1) [below of=A1] {$ \otimes^1 $};
  \node (B0) [left of=B1] {$\otimes^1 \circ \otimes^1$};
  \draw[->,  above right] (A0) to node {$ 1_{\otimes^1} $} (B1);  
  \draw[->,  left] (A0) to node {$\iota$} (B0);
  \draw[->,  below] (B0) to node {$ \gamma $} (B1);
 
\end{tikzpicture}\end{center}

\noindent Commutativity of this triangle in turn implies commutativity of the square above, and consequently of the square \textbf{E} in question. 

Finally if $T_{[a,sa]} = 1_I$, then all four vertices of the square in question are $1_I$, all four natural transformations are simply $1$, and commutativity follows.

\end{enumerate}

\vspace{.2cm}
We have been operating under the assumption that $\eta^{\bbb}_{\aaa}$ was of the form $\gamma \bullet \iiota$.  We must now consider when it is an identity.  First, let us suppose that $\eta^{\bbb}_{\aaa} = 1_{1_I}$. Thus $T_{\bbb} =  1_I$ and hence $A_{0,n} = I$.  It is a special fact about $\ctn$ that if $x_{0,n}=0$, every $1$-face of $x$ must also be $0$, and so $x$ is a degeneracy of the $1$-simplex $0$. This contradicts our assumption that $x$ is non-degenerate.  We need not worry about this case.

Now suppose that $\eta^{\bbb}_{\aaa} = 1_{1_A}$. This immediately implies that $T_{\bbb}=1_A$ and $T_{\aaa} = 1_A$. Consequently, the products comprising their domains -- $\prod_{a \in \aaa^-}A^{x_{a,sa}}$ and $\prod_{b\in\bbb^-}A^{x_{b,sb}}$ -- each contain only a single $A$, with the rest $I$. This then means that every product of functors $into$ those products must consist entirely of a product of $1_I$'s, except for the functor which targets the single $A$.   Hence, there exist unique $a' \in \aaa^-$ and unique $b'\in \bbb^-$ such that:

\[ T_{[a',sa']}  \neq 1_I \neq T_{[b',sb']} \]
\\
Figure $(\ddagger)$ then reduces to:

\begin{center}\begin{tikzpicture}[node distance=3.2cm]
\node (A1) {$\prod \limits_{i \in \n^-}A^{x_{i,i+1}}$};
\node (A2) [right of=A1] {$ $};
\node (A3) [right of=A2] {$A$};
\node (B1) [below of=A1] {$A$};
\node (B2) [right of=B1] {$ $};
\node (B3) [right of=B2] {$A$};

 \draw[->,  above] (A1) to node {$T_{[a',sa']}$} (A3);
 \draw[->,  left] (A1) to node {$T_{[b',sb']}$} (B1);
 \draw[->,  below left, bend left=25] (A1) to node {$ T =\otimes^{\sum \spi{x}}$} (B3);
 \draw[->,  below] (B1) to node {$ 1_A $} (B3);
 \draw[->,  right] (A3) to node {$ 1_A $} (B3);
 \draw[->,  below right, bend right=25] (B1) to node {$ 1_A $} (A3);

\end{tikzpicture}\end{center}

\noindent We must then verify:

\[\eta_{\aaa} \bullet \left(1 \circ \eta^{[a',sa']}_{\bbb \cap [a',sa']} \right) = \eta_{\bbb} \bullet \left( 1 \circ 1\right)\]
\\\\
Towards these end we note two further cases: $\eta^{[a',sa']}_{\bbb \cap [a',sa']} = \gamma \bullet \iiota$  or  $\eta^{[a',sa']}_{\bbb \cap [a',sa']} = 1_{T_{[a',sa']}}$.  In the first case, it must be that $T_{[a',sa']}$ is a tensor given by the number of $A$'s in its domain, hence, $T_{[a',sa']} = \otimes^{\sum \spi{x}}$.  We see then that $\eta^{[a',sa']}_{\bbb \cap [a',sa']} = \eta_{\bbb}$ and $\eta_{\aaa} = 1_{\otimes^{\sum \spi{x}}}$ as a result of the unitality axioms of the lax monoidal category.  Hence the equation above is verified.  In the second case, it must be that $T_{[a',sa']} = 1_A = T_{[b',sb']}$, and we see by inspection that $\eta_{\aaa} = \eta_{\bbb}$ once again verifying the equation.

\end{proof}

\medskip

\begin{proposition}

The assignment $(A,\otimes^n, \gamma, \iota) \mapsto \alpha: \ctn \to \hcn$ classifies lax monoidal categories.  That is, given a map $\alpha$, we can recover the data $(A,\otimes^n, \gamma, \iota)$.

\end{proposition}

\medskip

\begin{proof}

Given $\alpha$, we see that $\alpha_1(1)(\lbrace 0,1 \rbrace) = A^{1_{0,1}} = A^1 = A$, hence we have recovered the category $A$.  For $n\ge 2$, let $\mu \in \ctn_n$ be the $n$-simplex with $1$-faces $\mu_{p,q} =1$ for all $p,q \in \n$.  Then we have:

\[ \alpha_n(\mu)(\lbrace 0,n \rbrace \subset \n) = T^{\mu} = \otimes^{\sum \spi{\mu}} = \otimes^n\]

\medskip
\noindent Hence we have recovered the $n$-ary operations $\otimes^n$ for $n\ge 2$.  We can recover $\otimes^0$ from the non-degenerate $2$-simplex $u :0 \vee 0 \to 1 $.  We have:
\[ \alpha_2(u)(\lbrace 0,2 \rbrace \subset [2]) = T^u = \otimes^{\sum \spi{u}} = \otimes^0 \]
\medskip
\noindent  We can recover $\otimes^1$ from the non-degenerate $3$-simplex $l \in \ctn_3$ defined by the following $1$-faces:
\[ l_{0,1}=0 ,l_{1,2}=0, l_{2,3}=1, l_{0,2}=1, l_{1,3}=1, l_{0,3}= 1\] 
\medskip
\noindent
We have:
\[\alpha_3(l)(\lbrace 0,3 \rbrace \subset [3]) = T^{l}= \otimes^{\sum \spi{l}} = \otimes^1 \]
\medskip
This $3$ simplex $l$ will also recover the natural transformation $\iota$:
\[ \alpha_3(l)(\lbrace 0,3 \rbrace \subset \lbrace 0,1,3 \rbrace \subset [3]) = \eta^l_{\lbrace 0,1,3\rbrace} : T^{l}_{013} \circ \left(T^{l}_{01} \times T^{l}_{123}\right) \tto T^{l}_{[3]}\]
\medskip
\noindent As $l_{013} = l_{123} = s_1(1)$, $ T^{l_{013}} = T^{l_{123}} = 1_A$, and because $l_{01}=0$ is a degenerate $1$-simplex,   $T^{{l}_{01}}= 1_I$.  The definition of $\eta^l_{\lbrace 0,1,3\rbrace}$ becomes:

\[\eta^l_{\lbrace 0,1,3\rbrace} = \gamma_{1,2} \circ \left( \iota \times \iota \right): 1_A \circ 1_A \tto \otimes^1\]
\medskip
\noindent
By the unitality axiom for the lax monoidal category, this composite is simply $\iota$.

\medskip
Finally, as for recovering each $\gamma_{n,k_1,...,k_n}$, when $n\ge 2$ consider the simplex $x$ defined as follows.  Off its spine, every $1$-face $x_{p,q} = 1$.  Its spine, $\lbrace x_{i,i+1} \rbrace$, a finite sequence of $0$'s and $1$'s , will have two $0$'s for each $k_i = 0$, two $0$'s and a $1$ for each $k_i = 1$, and $k_i$ $1$'s for each $k_i \ge 2$, ordered with $k_1$'s digits first and $k_n$'s last.  Let $\hat{k_i} \in \mathbb{N}$ stand for the number of digits in the spine of $x$ associated with $k_i$, i.e $\hat{0} = 2$, $\hat{1} = 3$, $\hat{k_i}=k_i$ for $k_i\ge 2$.  Then $x \in \ctn_k$ with $k := \sum\limits_{1\le i \le n} \hat{k_i}$.  Writing $\hat{j} := \sum\limits_{1\le i \le j-1}\hat{k_i}$, we see that the functor $T^{x_{[\hat{j},\hat{j+1}]}}$ corresponding to the face $x_{[\hat{j},\hat{j+1}]}$ is precisely $\otimes^{k_j}$ by the above.  We have therefore: 
\[\alpha_{k}(\mu)\left(\lbrace 0, k \rbrace \subseteq \lbrace 0,\hat{k_1}, \hat{k_1}+\hat{k_2},..., k \rbrace \subseteq [k] \right) = \gamma_{n,k_1,...,k_n}: \otimes^n \circ (\otimes^{k_1} \times...\times \otimes^{k_n}) \tto \otimes^{k_1 + ... + k_n} \]

\medskip

\noindent The cases when $n=0$ or $n=1$ proceed similarly.

\end{proof}

\medskip
This concludes the proof of proposition \ref{laxclass}.
\medskip

\section{Classifying skew monoidal categories}

In this section, we will prove proposition \ref{skewclass}: that the nerve $\skn$ appearing in the classifcation result of \cite{third} embeds into $\hcn$.

\begin{definition}
The \textit{skew nerve} $\skn$ of the monoidal bicategory $\mathsf{Cat}$ is the simplicial set defined by the following:

\begin{enumerate}[\textbf{$\bullet$}]

\item There is a unique $0$-simplex, $*$.
\item $1$-simplices consist in categories $B_{01}$.
\item $2$-simplices consist in functors $B_{012} : B_{01}\times B_{12} \to B_{02}$
\item $3$-simplices are natural transformations $B_{0123}:B_{013}\circ (B_{012}\times 1)\tto B_{023} \circ (1 \times B_{123})$:

\begin{center}\begin{tikzpicture}[node distance=3cm]

  \node (A1) {$ B_{01} B_{12}B_{23}$} ;
  \node (A2) [right of=A1] {$B_{02} B_{23}$};
  \node (B1) [below of=A1] {$B_{01} B_{13}$};
  \node (B2) [below of=A2] {$B_{03}$};

  \draw[->, above] (A1) to node {$ B_{012}\times 1$} (A2);  
  \draw[->, left] (A1) to node {$ 1\times B_{123}$} (B1); 
 \draw[->, line width=.08cm, below right] (B1) to node {$ B_{0123}$} (A2); 
  \draw[-, white, line width=.06cm, below right] (B1) to node {$ $} (A2); 
  \draw[->, below ] (B1) to node {$ B_{013}$} (B2); 
  \draw[->, right ] (A2) to node {$ B_{023}$} (B2); 
  
\end{tikzpicture}\end{center}

\item A $4$-simplex consists in a quintuple of appropriately formed natural transformations making the following pentagon commute:

\begin{center}\begin{tikzpicture}[node distance=3cm]
 \node (A1) {$ B_{014} \circ (1\times B_{124}) \circ (1\times 1 \times B_{234}) $} ;
  \node (A2) [right of=A1] {$ $};
  \node (B2) [below of=A2] {$B_{024} \circ (B_{012} \times 1) \circ (1\times 1 \times B_{234}) $};
  \node (B1) [below of=A1] {$ $};
  \node (B0) [left of=B1] {$ B_{014} \circ (1 \times B_{134}) \circ (1 \times B_{123} \times 1)$};
  \node (C2) [below of=B2] {$B_{024} \circ (1 \times B_{234}) \circ (B_{012} \times 1 \times 1)$};
  \node (C1) [below of=B1] {$ $};
  \node (C0) [left of=C1] {$B_{034} \circ (B_{013} \times 1) \circ (1 \times B_{123} \times 1)$};
  \node (D1) [below of=C1] {$B_{034} \circ (B_{023} \times 1) \circ ( B_{012} \times 1 \times 1)$};

  \draw[-, line width=.08cm, right] (B2) to node {$ 1 $} (C2); 
  \draw[-, white, line width=.03cm, below right] (B2) to node {$ $} (C2);   
  
  \draw[->,  above right] (A1) to node {$B_{0124} \circ 1$} (B2);  
  \draw[->, below right] (C2) to node {$B_{0234}\circ 1$} (D1);
  \draw[->,  above left] (A1) to node {$1\circ(1 \times B_{1234})$} (B0);  
  \draw[->,  left] (B0) to node {$B_{0134}\circ 1$} (C0);
  \draw[->, below left] (C0) to node {$1\circ (B_{0123}\times 1)$} (D1);
 
\end{tikzpicture}\end{center} 

\item Higher-dimensional simplices are determined by $4$-coskeletality.

\end{enumerate}

\medskip

\end{definition}

\medskip

As in the case of the definition of a monoidal category, the pentagon law above gives rise to a coherence theorem which we now explain.  The $2$-faces of an $n$-simplex $B \in \skn$ consist in functors for every triple of numbers $p,q,r$ with $0 \le p < r < q \le n$: 

\[ B_{prq} : B_{pr} \times B_{rq} \to B_{pq}\]

\noindent Given a subset of indices $\ccc = \lbrace p = c_0 < ... < c_m = q \rbrace$, we then have a number of composite functors formed of these $2$ faces:
\[ \prod \limits_{c \in \ccc^-} B_{c,sc} \to B_{pq} \]

\noindent Moreover, we have a potential multitude of natural transformation $3$-faces mediating between such composites. For example, the pentagon diagram occurring in the definition of $\skn$ shows all of the composite $2$-face functors from $\prod_{i \in [4]^-} B_{i,i+1} \to B_{0,4}$, and shows all the $3$-face natural transformations between them.  

\medskip

In the context of a simplex $B \in \skn_n$, the coherence theorem says that there is at most one composite of $3$-face natural transformations between composites of $2$-face functors.  This is precisely the content of the commutativity of the pentagon above, and as in the theorem of \cite{second}, commutativity of pentagons gives the result in full generality.  Note also that every $2$-face functor composite $\prod_{c \in \ccc^-} B_{c,sc} \to B_{p,q}$ is the source of a composite of $3$-face transformations with target:
\[ m(B_{\ccc}):= B_{c_0 c_{m-1} c_m} \circ (B_{c_0 c_{m-2}c_{m-1}}\times 1) \circ ...\circ (B_{c_0 c_2 c_3} \times 1 )\circ ( B_{c_0 c_1 c_2} \times 1 ) \]

\medskip

\noindent  In a $4$-simplex $B\in \skn_4$, this $2$-face functor is just the bottom vertex of the pentagon above, $m(B)=B_{034} \circ (B_{023} \times 1) \circ ( B_{012} \times 1 \times 1)$.  Combined with the coherence theorem we have the following.  

\begin{remark}\label{mynoc}
Let $B \in \skn_n$ with $n\ge 2$, and $\ccc \subseteq \n$.  For each composite of $2$-face functors $T:\prod_{c\in \ccc^-} B_{c,sc} \to B_{p,q}$, there exists a unique natural transformation formed of composites of $3$-face transformations $T \tto m(B_{\ccc})$. 
\end{remark}

\medskip

Finally, we also have that $m(B) = m(s_i(B))$ for any $B \in \skn_n$ with $n \ge 2$ and $s_i: \skn_n \to \skn_{n+1}$ any degeneracy map.  Again recalling lemma \ref{degens}, we have:

\begin{align*}
m(s_i(B)) =& s_i(B)_{0,n,n+1 } \circ s_i(B)_{0,n-1,n } \circ ... \circ s_i(B)_{0,i,i+1 } \circ ... \circ s_i(B)_{0,1,2} \\
=& B_{\sigma_i(0,n,n+1)} \circ B_{\sigma_i(0,n-1,n)} \circ ... \circ s'(B_{0,i}) \circ ... \circ B_{\sigma_i(0,1,2)} \\
=& B_{0,n-1,n} \circ B_{0,n-2,n-1} \circ... \circ 1_{B_{0,i}} \circ ... \circ B_{0,1,2}\\
=& m(B)
\end{align*}

\begin{definition}\textbf{The Skew Nerve Embedding}

We will again rely on proposition \ref{backward}.  As we tend to write simplices in $\skn$ with capital letters $B$, we will recall that proposition in slighty different terminology.  To specify a map $\beta: \skn \to \hcn$, we must specify:

\begin{enumerate}
\item  A category for each $1$-simplex $B \in \skn_1$.  We will simply assign each $1$-simplex category $B$ to itself, as this satisfies $(1.a)$.

\item  A functor 

\[T(B): \prod \limits_{i \in \n^-} B_{i,i+1} \to B_{0,n}\]

\noindent For each $n \ge 1$ and $B \in \skn_n$.   Let $T(B) = 1_B$ when $n = 1$ in accordance with $(2.a)$, and otherwise let $T(B)=m(B)$ when $n \ge 2$, as by the above this is in accordance with $(2.b)$.  

\item  A natural transformation

\[ \eta^B_{\ccc}:  m(B_{\ccc}) \circ  \prod\limits_{c\in \ccc^-} m(B_{[c,sc]}) \tto m(B)\]

\noindent For each $n \ge 1$, $B \in \skn_n$ and $\lbrace 0,n \rbrace \subseteq \ccc \subseteq \n$.  For $1 \le n \le 2$, we take $\eta^B_{\ccc}:= 1_{B}$ in accordance with $(3.a)$.  Otherwise, take $\eta^B_{\ccc}$ to be the unique such transformation formed of composite of $3$-faces of $B$ in all cases given by remark \ref{mynoc}. It is easy to check that this assignment is in accordance with $(3.b)$ given that $m(s_i(B))=m(B)$.

\end{enumerate}

\end{definition}

This assignment gives rise to a map $\skn \to \hcn$ if equation \ref{key3} is satisfied for every $n \ge 4$ and $(\lbrace 0,n \rbrace \subset \aaa \subset \bbb \subset \n)$.  This is again the case by the uniqueness clause of \ref{mynoc}.

Therefore by \ref{backward} we have a map $\beta: \skn \to \hcn$ such that:

\begin{center} $\beta(B)(\lbrace 0,1 \rbrace) = B_{0,1}$, $\beta(B)(\{0,n\} \subseteq \n) = m(B)$, and $\beta(B)(\{0,n\} \subseteq \ccc \subseteq \n) = \eta^B_{\ccc}$. \end{center}

\begin{proposition}

The map $\beta_n:\skn_n \to \hcn_n$ is injective for each $n \ge 0$.  Hence $\beta$ is a faithful embedding.

\end{proposition}

\medskip

\begin{proof}

Let $B, D \in \skn_n$ and $\beta_n(B)=\beta_n(D)$.  For each subset $\lbrace p\le r\le q \rbrace \subset [n]$ we have $B_{pq} = \beta_n(B)(\lbrace p,q \rbrace)$ and $B_{prq} = m(B_{prq})= \beta_n(B)(\lbrace p,q \rbrace \subseteq \lbrace p,r,q \rbrace)$, hence $B$ and $D$ have precisely the same $1$ and $2$-faces.

Note that for $ \ccc = \lbrace c_0 < c_1 < c_2 < c_3 \rbrace$, we have that: \[ \beta_n(B)(\lbrace c_0,c_3 \rbrace \subset \ccc) = m(B_{\ccc})= B_{c_0c_2c_3} \circ (B_{c_0c_1c_2} \times 1) \] because $\beta$ commutes with face maps.  We therefore have:
\begin{align*}
\beta_n(B)(\lbrace c_0,c_3 \rbrace \subset \lbrace c_0,c_1,c_3 \rbrace \subset \ccc) &: m(B_{c_0c_1c_3}) \circ (m(B_{c_0c_1}) \times m(B_{c_1c_2c_3})) \tto m(B_{\ccc}) \\ 
= B_{c_0c_1c_2c_3} &: B_{c_0c_1c_3} \circ (1 \times B_{c_1c_2c_3}) \tto B_{c_0c_2c_3} \circ (B_{c_0c_1c_2} \times 1)
\end{align*}

\medskip\noindent This implies that $B$ and $D$ have the same $3$-faces.  As a result, commutative pentagons of such transformations occur in $B$ exactly when they occur in $D$, and so the two have the same $4$-faces.  They have the same $k$-faces for $k>4$ by $4-$coskeletality.
\end{proof}

\medskip
This concludes the proof of proposition \ref{skewclass}.
\medskip

\section{Classifying $\Sigma$-monoidal categories and general maps in $\mathsf{sSet}(\ctn,\hcn)$}

We have seen now that both lax monoidal categories (section $4$), skew monoidal categores (section $5$), and hence both monoidal and strict monoidal categories, all can be understood as maps in $\mathsf{sSet}(\ctn,\hcn)$.  The question we turn now to explore is then: what other monoidal-type categories do we find in this simplicial set?  We begin with a simple comparison of those $ \beta, \alpha,$ and $\phi:\ctn \to \hcn$ corresponding to skew monoidal categories, lax monoidal categories, and arbitrary maps respectively.  Recalling propositions \ref{forward} and \ref{backward}, we have that these three maps determine and are determined by three types of data:

\begin{enumerate}
\item A category $A^1=A$, where $A^0 = I$

\item A functor  
 \[T^x:\prod \limits_{i \in \n^-} A^{x_{i,i+1}} = A^{\sum \spi{x}} \to A\]

for each $n \ge 1$ and $x \in \ctn_n$ subject to stipulations $(2.a)$ and $(2.b)$ of proposition \ref{backward}.

\item A natural transformation
\[ \eta^x_{\ccc} : (T^{x_\ccc} \circ \prod\limits_{c \in \ccc^-}T^{x_{[c,sc]}}) \tto T^x \] 

for each face $n \ge 1$,  $x \in \ctn_n$, and $\lbrace 0,n \rbrace \subseteq \ccc \subseteq \n$ subject to $(3.a)$ and $(3.b)$ of proposition \ref{backward}.

\end{enumerate}

\noindent And, this data is equivalent to defining a map $\ctn \to \hcn$ whenever the natural transformations of item $(3)$ satisfy equation \ref{key3} of \ref{backward} for every $\lbrace 0,n \rbrace \subset \aaa \subset \bbb \subset \n$.  Let us then explore these three types of data associated to $\beta, \alpha$, and $\phi$.  Recall the five $2$-simplices of $\ctn$:  We have $s_0(0):0\vee 0 \to 0$, $s_0(1):0 \vee 1 \to 1$, $s_1(1):1 \vee 0 \to 1$, $u:0 \vee 0 \to 1$, and $m:1 \vee 1 \to 1$.  Here we write the simplex with $1$-faces $x_{0,1}$, $x_{1,2}$, and $x_{0,2}$ as $x_{0,1} \vee x_{1,2} \to x_{0,2}$.

\medskip

\begin{enumerate}

\item The category data associated to each of $\beta, \alpha$, and $\phi$ is just the monoidal-type category associated to the map.  

\medskip

\item The functor data associated to the skew monoidal category map $\beta$ is defined by a pair of functors, $\otimes^2:=T^{m}$, and $\otimes^0:=T^{u}$.  By the stipulations $(2.a)$ and $(2.b)$, this alone gives $T^x$ for every simplex $x \in \ctn_2$.  For an arbitrary $x \in \ctn_n$, $T^x$ is defined in terms of the image of its $2$-faces:
\[ T^x:= T^{x_{\lbrace 0,1,n \rbrace}} \circ ( 1 \times T^{x_{\lbrace 1,2,n \rbrace }}) \circ ... \circ (1 \times T^{ x_{\lbrace n-2,n-1,n\rbrace} })\]

\medskip

In contrast, the functor data associated to a lax monoidal category map $\alpha$ is defined by functors $\otimes^n:=T^{\mu}$ for every $n \ge 0$ where $\mu \in \ctn_n$ is given by $\mu_{p,q}=1$ for all $0 \le p<q \le n$.  For arbitrary $x \in \ctn$, $T^x:= \otimes^{\sum \spi{x}}$.  

\medskip  Finally, $\phi$ may specify apriori unrelated functors $T^x:A^{\sum \spi{x}} \to A$ for every non-degenerate $x$.  Note that there are a countable infinite number of non-degenerate simplices $x$ with $\sum\spi{x} = n$ for every $n \ge 0$,\footnote{ Take a simplex with $x_{i,i+1} = 1$ for $i \le n-1$, $x_{i,i+1}=0$ for $i \ge n$, and $x_{p,q}=1$ for $q-p\ge 2$.  There is one simplex of this form in every dimension $\ge n$, it is non-degenerate, and $\sum \spi{x} = n$.} and thus, $\phi$ may include the data of an arbitrary (countable) number of $n$-ary functors for each $n$.

\medskip

\item Key differences in the natural transformation data associated to $\beta, \alpha,$ and $\phi$, can already be seen when considering only $3$-simplices $x\in\ctn_3$.  Let us represent such a $3$-simplex with $1$-faces $x_{p,q}$, $0 \le p<q \le 3$ as the diagram:

\begin{center} \begin{tikzpicture}[node distance=1.25cm]
  \node (A1) {$ x_{0,1} \vee x_{1,2} \vee x_{2,3}$} ;
  \node (A2) [right of=A1] {$ $};
  \node (A3) [right of=A2] {$ x_{0,1} \vee x_{1,3}$};
  \node (B2) [below of=A2] {$x$};
  \node (B1) [below of=A1] {$ $};
  \node (B3) [below of=A3]{$ $};
  \node (C1) [below of=B1] {$x_{0,2} \vee x_{2,3}$};
  \node (C2) [below of=B2] {$ $};
  \node (C3) [below of=B3] {$x_{0,3}$};
  
  \draw[->,  above] (A1) to node {$ $} (A3);  
  \draw[->,  left] (A1) to node {$ $} (C1);
  \draw[->, below] (C1) to node {$ $} (C3);
  \draw[->, right] (A3) to node {$ $} (C3);
\end{tikzpicture}\end{center}

\medskip

 \noindent Consider then the pair of simplices $l$ and $r \in \ctn_3$:

\medskip

\begin{center} \begin{tikzpicture}[node distance=1.25cm]
  \node (A1) {$ 0 \vee 0 \vee 1$} ;
  \node (A2) [right of=A1] {$ $};
  \node (A3) [right of=A2] {$ 0 \vee 1$};
  \node (B2) [below of=A2] {$l$};
  \node (B1) [below of=A1] {$ $};
  \node (B3) [below of=A3]{$ $};
  \node (C1) [below of=B1] {$1 \vee 1$};
  \node (C2) [below of=B2] {$ $};
  \node (C3) [below of=B3] {$1$};
  
  \draw[->,  above] (A1) to node {$ $} (A3);  
  \draw[->,  left] (A1) to node {$ $} (C1);
  \draw[->, below] (C1) to node {$ $} (C3);
  \draw[->, right] (A3) to node {$ $} (C3);

\end{tikzpicture}\hspace{1cm}\begin{tikzpicture}[node distance=1.25cm]
  \node (A1) {$ 1 \vee 0 \vee 0$} ;
  \node (A2) [right of=A1] {$ $};
  \node (A3) [right of=A2] {$ 1 \vee 1$};
  \node (B2) [below of=A2] {$r$};
  \node (B1) [below of=A1] {$ $};
  \node (B3) [below of=A3]{$ $};
  \node (C1) [below of=B1] {$1 \vee 0$};
  \node (C2) [below of=B2] {$ $};
  \node (C3) [below of=B3] {$1$};
  
  \draw[->,  above] (A1) to node {$ $} (A3);  
  \draw[->,  left] (A1) to node {$ $} (C1);
  \draw[->, below] (C1) to node {$ $} (C3);
  \draw[->, right] (A3) to node {$ $} (C3);

\end{tikzpicture}\end{center}

\noindent Associated to each of $\beta, \alpha,$ and $\phi$ is a functor $T^l:A^{\sum \spi{l}} = IIA=A \to A$, as well as functors $T^{l_{\lbrace 0,1,2 \rbrace}}=T^u$, $T^{l_{\lbrace 0,1,3 \rbrace}}=T^{s_0(1)}=1_A$, $T^{l_{\lbrace 0,2,3 \rbrace}} = T^m$, and $T^{l_{\lbrace 1,2,3 \rbrace}}= T^{s_0(1)}=1_A$ corresponding to the four arrows making up the edges of the above squares.   There are also the pair of natural transformations: 

\[ \eta^l_{\lbrace 0,2,3 \rbrace}: T^{l_{\lbrace 0,2,3\rbrace}} \circ (T^{l_{\lbrace 0,1,2 \rbrace}} \times T^{l_{\lbrace 2,3 \rbrace}}) \tto T^{l}\]
\[ \eta^l_{\lbrace 0,1,3 \rbrace}: T^{l_{\lbrace 0,1,3\rbrace}} \circ (T^{l_{\lbrace 0,1 \rbrace}} \times T^{l_{\lbrace 1,2,3 \rbrace}}  ) \tto T^{l}\]

\medskip

\noindent We get similar data for $r$, and we can represent all of it succinctly in the following two diagrams:

\begin{center} \begin{tikzpicture}[node distance=2cm]
  \node (A1) {$ IIA$} ;
  \node (A2) [right of=A1] {$ $};
  \node (A3) [right of=A2] {$ IA$};
  \node (B2) [below of=A2] {$ $};
  \node (B1) [below of=A1] {$ $};
  \node (B3) [below of=A3]{$ $};
  \node (C1) [below of=B1] {$AA$};
  \node (C2) [below of=B2] {$ $};
  \node (C3) [below of=B3] {$A$};
  
  \draw[->,  above] (A1) to node {$1_I \times T^{l_{\lbrace 1,2,3 \rbrace}}$} (A3);  
  \draw[->,  left] (A1) to node {$T^{l_{\lbrace 0,1,2 \rbrace}} \times 1_A$} (C1);
  \draw[->, below] (C1) to node {$T^{l_{\lbrace 0,2,3 \rbrace}}$} (C3);
  \draw[->, right] (A3) to node {$T^{l_{\lbrace 0,1,3 \rbrace}}$} (C3);
  \draw[->, above left] (A1) to node {$T^l$} (C3);
  \draw[->, below left] (A3) to node {$ $} (B2);
 \draw[->, line width=.08cm, above left] (A3) to node {$   \eta^l_{\lbrace 0,1,3 \rbrace}$} (B2); 
 \draw[-, white, line width=.03cm, below right] (A3) to node {$ $} (B2);   
 \draw[->, line width=.08cm, below right] (C1) to node {$  \eta^l_{\lbrace 0,2,3 \rbrace}$} (B2); 
 \draw[-, white, line width=.03cm, below right] (C1) to node {$ $} (B2);   

\end{tikzpicture}\begin{tikzpicture}[node distance=2cm]
  \node (A1) {$ AII$} ;
  \node (A2) [right of=A1] {$ $};
  \node (A3) [right of=A2] {$ AA$};
  \node (B2) [below of=A2] {$ $};
  \node (B1) [below of=A1] {$ $};
  \node (B3) [below of=A3]{$ $};
  \node (C1) [below of=B1] {$AI$};
  \node (C2) [below of=B2] {$ $};
  \node (C3) [below of=B3] {$A$};
  
  \draw[->,  above] (A1) to node {$1_A \times T^{r_{\lbrace 1,2,3 \rbrace}}$} (A3);  
  \draw[->,  left] (A1) to node {$T^{r_{\lbrace 0,1,2 \rbrace}} \times 1_I$} (C1);
  \draw[->, below] (C1) to node {$T^{r_{\lbrace 0,2,3 \rbrace}}$} (C3);
  \draw[->, right] (A3) to node {$T^{r_{\lbrace 0,1,3 \rbrace}}$} (C3);
  \draw[->, above left] (A1) to node {$T^r$} (C3);
  \draw[->, below left] (A3) to node {$ $} (B2);
 \draw[->, line width=.08cm, above left] (A3) to node {$  \eta^r_{\lbrace 0,1,3 \rbrace}$} (B2); 
 \draw[-, white, line width=.03cm, below right] (A3) to node {$ $} (B2);   
 \draw[->, line width=.08cm, below right] (C1) to node {$  \eta^r_{\lbrace 0,2,3 \rbrace}$} (B2); 
 \draw[-, white, line width=.03cm, below right] (C1) to node {$ $} (B2);   

\end{tikzpicture}\end{center}

The remaining functors associated to the skew monoidal category map $\beta$ are defined $T^l:=1_A$ and $T^r:= \otimes^2 \circ (1_A \times \otimes^0)$.  Naming $\lambda := \eta^l_{0,2,3}$ and $\rho:= \eta^r_{0,2,3}$, we get the following diagrams:

\begin{center} \begin{tikzpicture}[node distance=2cm]
  \node (A1) {$ IIA$} ;
  \node (A2) [right of=A1] {$ $};
  \node (A3) [right of=A2] {$ IA$};
  \node (B2) [below of=A2] {$ $};
  \node (B1) [below of=A1] {$ $};
  \node (B3) [below of=A3]{$ $};
  \node (C1) [below of=B1] {$AA$};
  \node (C2) [below of=B2] {$ $};
  \node (C3) [below of=B3] {$A$};
  
  \draw[->,  above] (A1) to node {$1_I \times 1_A$} (A3);  
  \draw[->,  left] (A1) to node {$ \otimes^0 \times 1_A$} (C1);
  \draw[->, below] (C1) to node {$ \otimes^2 $} (C3);
  \draw[->, right] (A3) to node {$1_A$} (C3);
  \draw[->, above] (A1) to node {$ $} (C3);
  \draw[->, below left] (A3) to node {$ $} (B2);
 \draw[=, line width=.08cm, above left] (A3) to node {$ 1 $} (B2); 
 \draw[-, white, line width=.03cm, below right] (A3) to node {$ $} (B2);   
 \draw[->, line width=.08cm, below right] (C1) to node {$ \lambda $} (B2); 
 \draw[-, white, line width=.03cm, below right] (C1) to node {$ $} (B2);   

\end{tikzpicture}\hspace{1cm}\begin{tikzpicture}[node distance=2cm]
  \node (A1) {$ AII$} ;
  \node (A2) [right of=A1] {$ $};
  \node (A3) [right of=A2] {$ AA$};
  \node (B2) [below of=A2] {$ $};
  \node (B1) [below of=A1] {$ $};
  \node (B3) [below of=A3]{$ $};
  \node (C1) [below of=B1] {$AI$};
  \node (C2) [below of=B2] {$ $};
  \node (C3) [below of=B3] {$A$};
  
  \draw[->,  above] (A1) to node {$1_A \times \otimes^0$} (A3);  
  \draw[->,  left] (A1) to node {$1_A \times 1_I$} (C1);
  \draw[->, below] (C1) to node {$1_A$} (C3);
  \draw[->, right] (A3) to node {$\otimes^2$} (C3);
  \draw[->, above] (A1) to node {$ $} (C3);
  \draw[->, below left] (A3) to node {$ $} (B2);
 \draw [=, line width=.08cm, above left] (A3) to node {$ 1 $} (B2); 
 \draw[-, white, line width=.03cm, below right] (A3) to node {$ $} (B2);   
 \draw[->, line width=.08cm, below right] (C1) to node {$ \rho $} (B2); 
 \draw[-, white, line width=.03cm, below right] (C1) to node {$ $} (B2);   

\end{tikzpicture}\end{center}

\noindent  We see then $\lambda$ and $\rho$ arise naturally in this context.  Associated to the lax monoidal category map $\alpha$ are instead $T^l = \otimes^1 = T^r$ and the following diagrams:

\begin{center} \begin{tikzpicture}[node distance=2cm]
  \node (A1) {$ IIA$} ;
  \node (A2) [right of=A1] {$ $};
  \node (A3) [right of=A2] {$ IA$};
  \node (B2) [below of=A2] {$ $};
  \node (B1) [below of=A1] {$ $};
  \node (B3) [below of=A3]{$ $};
  \node (C1) [below of=B1] {$AA$};
  \node (C2) [below of=B2] {$ $};
  \node (C3) [below of=B3] {$A$};
  
  \draw[->,  above] (A1) to node {$1_I \times 1_A$} (A3);  
  \draw[->,  left] (A1) to node {$\otimes^0 \times 1_A$} (C1);
  \draw[->, below] (C1) to node {$\otimes^2$} (C3);
  \draw[->, right] (A3) to node {$1_A$} (C3);
  \draw[->, above left] (A1) to node {$\otimes^1$} (C3);
  \draw[->, below left] (A3) to node {$ $} (B2);
 \draw[->, line width=.08cm, above left] (A3) to node {$  \iota $} (B2); 
 \draw[-, white, line width=.03cm, below right] (A3) to node {$ $} (B2);   
 \draw[->, line width=.08cm, below right] (C1) to node {$ \gamma_{2,0,1}$} (B2); 
 \draw[-, white, line width=.03cm, below right] (C1) to node {$ $} (B2);   

\end{tikzpicture}\hspace{1cm}\begin{tikzpicture}[node distance=2cm]
  \node (A1) {$ AII$} ;
  \node (A2) [right of=A1] {$ $};
  \node (A3) [right of=A2] {$ AA$};
  \node (B2) [below of=A2] {$ $};
  \node (B1) [below of=A1] {$ $};
  \node (B3) [below of=A3]{$ $};
  \node (C1) [below of=B1] {$AI$};
  \node (C2) [below of=B2] {$ $};
  \node (C3) [below of=B3] {$A$};
  
  \draw[->,  above] (A1) to node {$1_A \times \otimes^0 $} (A3);  
  \draw[->,  left] (A1) to node {$1_A \times 1_I$} (C1);
  \draw[->, below] (C1) to node {$1_A$} (C3);
  \draw[->, right] (A3) to node {$\otimes^2$} (C3);
  \draw[->, above left] (A1) to node {$\otimes^1$} (C3);
  \draw[->, below left] (A3) to node {$ $} (B2);
 \draw[->, line width=.08cm, above left] (A3) to node {$ \gamma_{2,1,0}$} (B2); 
 \draw[-, white, line width=.03cm, below right] (A3) to node {$ $} (B2);   
 \draw[->, line width=.08cm, below right] (C1) to node {$\iota$} (B2); 
 \draw[-, white, line width=.03cm, below right] (C1) to node {$ $} (B2);   

\end{tikzpicture}\end{center}

\medskip

\noindent So we see the transformation $\iota$ and two of the many $\gamma$ transformations arise naturally in this context as well.  Finally, for the general map $\phi$, there are, apriori no relationships between these various functors and natural transformations, aside from those relationships contributed by the other data of $\phi$.  

\end{enumerate}

\medskip

We have seen in particular that arbitrary maps $\phi$ may include the data of up to a countable number of distinct $n$-ary functors for each $n$.  Such a monoidal-type category is reminiscent of the following definition of \cite{fifth} which we present only informally here:

\begin{definition}

A $\Sigma$-\textit{monoidal category} $(A, \Sigma, \gamma)$ consists in a category $A$, a countable\footnote{The definition in \cite{fifth} considers aribtrary cardinality.} set $\Sigma_k$ of $k$-ary functors $A^k \to A$ for each $k \ge 0$, and natural isomorphisms $\gamma$ between each possible composite of functors of the same total arity. Each composition of these natural isomorphisms with the same domain and codomain must be equal.

\end{definition}

There is no canonical way of associating a map $\sigma \in\mathsf{sSet}(\ctn,\hcn)$ to a $\Sigma$-monoidal category $(A, \Sigma, \gamma)$ because the definition presents all $k$-ary functors in $\Sigma_k$ as morally indistinguishable, whereas the $k$-ary functors arising in the image of a map $\ctn \to \hcn$ can be distinguished in many ways, for example, by the dimension of the simplex $x$ mapping to that functor.  However, given a surjective function 
\[h_k: \lbrace x~ | ~ x \text{ is nondegenerate with dimension } \ge 2 \text{ and } \sum \spi{x}=k \rbrace \to \Sigma_k\]

\noindent for each $k\ge 0$, we can define $\sigma: \ctn \to \hcn$ by proposition \ref{backward}.

\begin{definition}\textbf{The $\Sigma$-Monoidal Classifying Map}

\begin{enumerate}

\item Let $A^x:= A$ if $x=1$ and $I$ otherwise.
\item Let $T^{x}:= h_{\sum \spi{x}}(x): A^{\sum \spi{x}} \to A$ if $x$ is non-degenerate with dimension $\ge 2$.  Otherwise $T^x$ is defined in accordance with stipulations $(2.a)$ and $(2.b)$.
\item Let $\eta^x_{\ccc}: T^{x_{\ccc}} \circ \prod_{c \in \ccc^-} T^{x_{[c,sc]}} \tto T^x$ be the unique $\gamma$ natural isomorphism guaranteed by the definition of $\Sigma$-monoidal category if $n \ge 3$, and $1_{T^x}$ otherwise, as this is in accordance with $(3.a)$ and $(3.b)$.
\end{enumerate}

\end{definition}

\medskip

\noindent The commutativity of equation \ref{key3} is implied directly by the commutativity of the isomorphisms $\gamma$ given in the definition, and hence this data gives rise to a map $\sigma:\ctn \to \hcn$.  

\medskip

This assignment also classifies $\Sigma$-monoidal categories.

\begin{proposition}\label{classsigma}

The assignment $(A,\Sigma, \gamma) \mapsto \sigma: \ctn \to \hcn$ classifies $\Sigma$-monoidal categories.  That is, given the map $\sigma$, we can recover the data $(A, \Sigma, \gamma)$.

\end{proposition}

We get $A$ from item (1), and each map from each $\Sigma_k$ from item (2), using the fact that $h_k$ is assumed to be surjective.  We must therefore only show that \textit{every} natural isomorphism implied by the definition of $\Sigma$-monoidal categories can be generated -- via horizontal composition, vertical composition, and by product -- by the transformations $\eta^x_{\ccc}$ of item (3). We need a few lemmas which will incidently reveal some additional structure of general maps $\phi$.

\begin{lemma}

Let  $\phi: \ctn \to \hcn$ be a general map.  Suppose $x \in \ctn_m$ such that $\sum \spi{x} = 1$.  Then there is a transformation $E:1_A \tto T^x$ generated by transformations $\eta^y_{\ccc}$ in the image of $\phi$. 

\end{lemma}

\medskip

\begin{proof}

Let $x \in \ctn_m$ with $\sum \spi{x} = 1$.  Let $x_{i-1,i}=1$ be the unique such $1$-face in $\spi{x}$.  If $i \neq 1$, then:
\[ \eta^x_{0,1,...,i-2,i,i+1,...,m}: T^{x_{0,1,...,i-2,i,i+1,...,m}} \circ (1_I \times ... \times 1_I \times 1_A \times 1_I \times ... \times 1_I) \tto T^x \]
\medskip

\noindent Here, we get $T^{x_{i-2,i-1,i}}=1_A$ because $x_{i-2,i-1,i}$ is necessarily $s_0(1)$, and hence degenerate.  So we get a transformation $T^{x_{0,1,...,i-2,i,i+1,...,m}} \tto T^x$.  On the other hand if $i=1$ we have:
\[ \eta^x_{0,2,...,m}: T^{x_{0,2,...,m}} \circ (1_A \times 1_I \times ... \times 1_I) \tto T^x\]
\medskip

\noindent This time $T^{x_{0,1,2}} = 1_A$ because $x_{0,1,2} = s_1(1)$.  In each case, the face $x_{0,1,...,i-2,i,i+1,...,m}$ and $x_{0,2,...,m}$ both have a single $1$ on the spine, and one fewer $0$.  By repeating the argument on these faces inductively, we get a composite:

\[ E: T^{x_{0,m}}=1_A \tto T^x \]

\end{proof}

\medskip

\begin{lemma}

Let  $\phi:\ctn \to \hcn$ be a general map.  Let $x \in \ctn_m$ with $\sum \spi{x}= 0$ and $x$ not a degeneracy of $0$.  Then there is a transformation $E: T^u \tto T^{x}$ generated by transformations $\eta^y_{\ccc}$ in the image of $\phi$.

\end{lemma}

\medskip

\begin{proof}

We consider two cases.  If $x_{0,1,2} = s_0(0)$ then $T^{x_{0,2,3,...,n}}=T^x$, and we are finished by induction.  If on the other hand $x_{0,1,2}=u$, then $\sum \spi{x_{0,2,3,...,n}} = 1$ and by the previous lemma, we have a transformation $F: 1_A \tto T^{x_{0,2,3,...,n}}$.  We have:

\[ \eta^x_{0,2,3,...,n}:  T^{x_{0,2,3,...,n}} \circ (T^u \times 1_I \times ... \times 1_I) \tto T^x\]
\[ E:= \eta^x_{0,2,3,...,n} \bullet (E \circ 1 ): T^u \tto T^x \]

\end{proof}

\medskip

\begin{lemma}

Let  $\phi:\ctn \to \hcn$ be a general map.  For $n\ge 2$, let $\mu \in \ctn_n$ be given by $\mu_{p,q}=1$ for all $0 \le p < q \le n$, and let $T^{\mu}$ be the associated $n$-ary functor.  Let $x \in \ctn_m$ with $\sum \spi{x}=n$.  Then there is a transformation $E: T^x \tto T^{\mu}$ generated by transformations $\eta^y_{\ccc}$ in the image of $\phi$.

\end{lemma}

\medskip

\begin{proof}

Now, let $\mu \in \ctn_n$ as above, and let $x \in \ctn_m$ with $m\ge n \ge 2$ with $\sum \spi{x} = n$.  Suppose that $x_{i-1,i}=1$ if and only if $i \in I = \lbrace i_1,...,i_n\rbrace$.  Then the faces $\sum \spi{x_{0,...,i_1}} =\sum \spi{x_{i_1,...,i_2,}} = ... = \sum \spi{x_{i_{n-1},...,i_n-1,i_n,...,m}} = 1$ . Thus we have:

\[ E_1 \times ... \times E_n : 1 \tto  T^{x_{0,...,i_1}}\times ... \times T^{x_{i_{n-1},...,m}}\]

\medskip

We also have that $x_{0,i_1,...,i_{n-1},m} = \mu$ because $\mu$ is the unique $n$-simplex with $\sum \spi{\mu} = n$.  We have then:

\[ \eta^x_{0,i_1,i_2,...,i_{n-1},m} : T^{x_{0,i_1,...,i_{n-1},m}} \circ \left(T^{x_{0,...,i_1}}\times ... \times T^{x_{i_{n-1},...,m}}\right) \tto T^x\]  

\[  E:= \eta^x_{0,i_1,i_2,...,i_{n-1},m} \bullet (1 \circ E_1 \times ... \times E_n): (T^{\mu} \circ 1) \tto (T^{\mu} \circ  T^{x_{0,...,i_1}}\times ... \times T^{x_{i_{n-1},...,m}}) \tto T^x\] 

\end{proof}

\medskip

\begin{proof}(of proposition \ref{classsigma})
Given a map $\sigma$ in the image of the assignment, we must show that every natural isomorphism implied by the definition of $\Sigma$-monoidal categories can be generated by transformations $\eta^x_{\ccc}$ in the image of $\sigma$.  As $\sigma$ assigns to each $\eta^x_{\ccc}$  a natural isomorphism, the above three lemmas imply that for every $k \ge 0$, each pair of elements $f,f' \in \Sigma_k$ are isomorophic to one another via isomorphisms in the image of $\sigma$.   Let $\mu^n \in \ctn_n$ denote the unique $n$-simplex with $\sum \spi{x}=n$, i.e $\mu^n_{p,q}=1$ for all $0 \le p< q \le n$.  Now, given $n = n_1 + ... + n_k$, $f_i \in \Sigma_{n_i}, f \in \Sigma_k$, and $g \in \Sigma_n$, we must show that the isomorphism $\gamma: f \circ (f_1 \times ... \times f_k) \tto g$ can be generated by $\eta^x_{\ccc}$'s.  Without loss of generality we can assume that $f_i  = T^{\mu^{n_i}}$ if $n_i \ge 2$, $f_i = 1_A$ if $n_i = 1$, and $f_i = T^u$ if $n_i = 0$ and similarly for $f$ and $g$.  Let $x$ be the simplex with two consecutive $0$'s on its spine for each $n_i = 0$, with $n_i$ consecutive $1$'s on its spine for each $n_i \ge 1$, and all other $1$-faces $x_{p,q} = 1$.   Let the total dimension of $x$ be $m$.  Then choosing $\ccc$ with $\lbrace0,m \rbrace \subseteq \ccc \subseteq [m]$ so that consecutive indicies of $\ccc$ correspond to the part of $\spi{x}$ corresponding to each $n_i$, we get the transformation:

\[ \eta^x_{\ccc} : f \circ (f_{n_1} \times ... \times f_{n_k}) \tto T^x\] 

Note that $x_{\ccc} = \mu^k$ as it has only $1$'s on its spine, and so $T^{x_{\ccc}} = f$.  Finally, as $\sum \spi{x} = n$, we get a transformation $E: T^x \tto g = T^{\mu^n}$.  Together these show the isomorphism between these $n$-ary functors to be generated by $\sigma$.

\end{proof}

\medskip

This concludes the proof of propostion \ref{sigmaclass}

\medskip

We might therefore think of general maps $\phi:\ctn \to \hcn$ as specifying a category along with an arbitrary (countable) number of $n$-ary functors for every $n$ which are not necessarily isomorophic.  They are, however, related by an intricate web of natural transformations which mirror aspects of the structure of $\ctn$. Perhaps one might therefore think of the data associated to an abitrary map $\ctn \to \hcn$ along with the equation \ref{key3} as indicating the necessary structure and coherence needed to go about weakening the definition of $\Sigma$-monoidal category to not require natural isomorphisms between all functors with the same total arity.

\bibliography{laxskewclass13}
\bibliographystyle{unsrt}

\end{document}